\documentclass[a4paper]{article}
\usepackage{amscd,amsfonts,amsmath,amssymb,amsthm}
\usepackage[english]{babel}
\usepackage{bbold}
\usepackage{caption}
\usepackage{color}
\usepackage{dsfont}
\usepackage{enumerate}
\usepackage{enumitem}
\usepackage{hyperref}
\usepackage{graphics}
\usepackage{graphicx}
\usepackage{mathrsfs}
\usepackage{mathtools}
\usepackage{mathabx}        
\mathtoolsset{showonlyrefs} 
\usepackage{pxfonts}
\usepackage{times}
\usepackage{xspace}

\usepackage[left=20mm,right=20mm, top=3cm, bottom=3cm]{geometry}

\newcommand{\Z}{\mathbb{Z}}
\newcommand{\R}{\mathbb{R}}
\renewcommand{\d}{\ \mathrm{d}}

\newcommand{\curl}{\mathrm{curl}}

\renewcommand{\d}{{\,\mathrm{d}}}

\newcommand{\domain}{D}

\newcommand{\macrogridsize}{H}

\newcommand{\Etensor}{\mathbf{C}}
\newcommand{\displaceMacro}{U}
\newcommand{\displaceMicro}{u}
\newcommand{\pf}{v}
\newcommand{\energy}{E}
\newcommand{\Eelast}{\energy_\text{elast}}
\newcommand{\Epot}{\energy_\text{pot}}
\newcommand{\Etot}{\energy_\text{tot}}
\newcommand{\Volume}{\text{Vol}}
\newcommand{\Perimeter}{\text{Per}}
\newcommand{\DoubleWell}{W}
\newcommand{\gridSizeSplineChart}{\tau}
\newcommand{\splineFESpace}{W_{\gridSizeSplineChart}}
\DeclareMathOperator*{\argmin}{arg\,min}

\theoremstyle{definition}

\theoremstyle{remark}

\theoremstyle{plain}

\makeatletter
\DeclareRobustCommand\onedot{\futurelet\@let@token\@onedot}
\def\@onedot{\ifx\@let@token.\else.\null\fi\xspace}
\def\eg{\emph{e.g}\onedot} 
\def\ie{\emph{i.e}\onedot} 
\def\cf{\emph{cf}\onedot}  

\def\wrt{w.r.t\onedot}     
        \def\st{s.t\onedot}
\def\etal{\emph{et al}\onedot}

\makeatother

\title{Two-scale elastic shape optimization for additive manufacturing}

\author{Sergio Conti$^1$, Martin Rumpf$^{1,2}$, Stefan Simon$^2$\\[3mm]
\small $^1$
    Institute for Applied Mathematics, University of Bonn, Endenicher Allee 60, 53115 Bonn, Germany\\
    \small
    $^2$    Institute for Numerical Simulation, University of Bonn, Endenicher Allee 60, 53115 Bonn, Germany
}
\date{\today}

\begin{document}

\maketitle

\begin{abstract}
In this paper, a two-scale approach for elastic shape optimization of fine-scale structures in additive manufacturing is investigated.
To this end, a free material optimization is performed on the macro-scale
restricting to isotropic
elasticity tensors in a set of microscopically realizable tensors.
A database of these realizable tensors and their cost values 
is obtained with a shape and topology optimization procedure on microscopic cells, working within a fixed set of elasticity tensors samples.
This microscopic optimization takes into account manufacturability constraints via predefined material bridges to neighbouring cells at the faces of the microscopic fundamental cell. 
For the actual additive manufacturing on a chosen fine scale, a piecewise constant elasticity tensor ansatz on the grid cells of a macroscopic mesh is applied. 
The macroscopic optimization is performed in an efficient online phase, during which the associated cell-wise optimal material patterns are retrieved from the database that was computed offline.
For that, the set of admissible realizable elasticity tensors is parametrized using tensor product cubic B-splines over the unit square matching the precomputed samples. 
This representation is then efficiently used in an interior point method for the free material optimization on the macro-scale.
\end{abstract}

\section{Introduction} \label{sec:intro}

We study shape optimization in linear elasticity 
with the constraint that the resulting shape should be additively manufactured.
To this end, we combine a free material optimization approach on the macro-scale with a topology and shape optimization approach on the micro-scale.  

The macroscopic problem selects in each discretization cell an isotropic elasticity tensor within an admissible set
determined by realizable microstructures.
Besides a macroscopic elasticity problem, for the local cost of an elasticity tensor, we take into account the volume occupied by hard material and the perimeter as a regularizer, which can be interpreted as introducing manufacturing costs for the production and processing of the object surface.
Note that the volume and the perimeter arise from the corresponding microcell.
Thus, in order to efficiently solve the macroscopic problem, the information of the microscopic problems
are stored in a precomputed database.

The microscopic problem determines the set of elasticity tensors that can be realized as homogenized tensors of a specific class of microscopic material patterns. 
This class is tailored to manufacturability and in particular ensures compatibility between adjacent macroscopic cells by requiring the presence of predefined material bridges at all boundaries. 
We stress that in standard two-scale approaches based on the theory of homogenization microstructures are only required to be periodic, so that mixing different microstructures leads to discontinuities at the boundaries, unless boundary layers on an additional, intermediate scale are introduced. 
Both options hamper manufacturability and strongly limit the applicability of the obtained patterns. 
Note that also in the case of two elastic phases with different, non degenerate material properties, a simple connection of adjacent cells without appropriate boundary conditions does not allow a proper transfer of stresses.
Because of the predefined material bridges, our approach, in contrast, ensures that any two microstructures are compatible on the boundaries of the microcells.
This directly allows manufacturability without post-processing of the original microstructures.

Numerically, the microscopic problem is very time consuming, and therefore it is performed offline. 
It consists of a shape and topology optimization procedure in which the effective elasticity tensors and the associated local costs are precomputed. 
They are then used for efficient online free material optimization on the macro-scale,
where, for simplicity, we restrict to locally isotropic elasticity tensors.

In general, for elastic bodies subjected to external surface loadings, 
it is well-known that shape optimization problems typically lead to the formation of microstructures, see \eg the textbook by Allaire~\cite{Al02}. 
More precisely, the associated variational problem is in general ill-posed. 
As a remedy, relaxation theory is usually applied with composite materials on a micro-scale and effective macroscopic elasticity tensors resulting from a microscopic mixture of the involved constituents.
Homogenization theory as presented in Cioranescu and Donato~\cite{DoCi99} or Milton~\cite{Mi02} is then applied to derive the effective material properties of a composite elastic material with two phases, \eg soft and hard.
Using the G-closure theory as detailed in Cherkaev~\cite{Ch00a}, it can be shown that a nested lamination structure realizes the minimal value of the compliance functional. 
These optimal constructions are usually not unique and entirely different microstructures might lead to identical macroscopic material properties.
Most of these structures are purely theoretical and can not be realized via additive manufacturing.

\paragraph*{Additive manufacturing.}
Driven by the rapid development of 3D printing devices, 
the design of microscopic material patterns, 
which can be used in an additive manufacturing approach as building blocks of a macroscopic workpiece,
recently attracted significant attention.
Panetta~\etal~\cite{PaZhMa15} proposed a class of 3D printable fine-scale structures on cubic fundamental cells 
consisting of a topologically restricted truss geometry with trusses of varying thicknesses and flexible node positions.
Thereby a tetrahedral symmetry was taken into account on the fundamental cells. This does not necessarily lead to isotropic homogenized elasticity tensors as mentioned in the paper.
Schumacher \etal~\cite{ScBiRy15} used harmonic displacements introduced in \cite{KhMuOw09} to optimize the cell microstructures.
Martinez~\etal~\cite{MaDuLe16} proposed so called Voronoi foams to generate printable 
microstructures.
In \cite{ZhSkCh17}, Zhu~\etal solved a macroscopic material optimization problem making use of a level set representation of microscopic shapes with effective macroscopically isotropic material parameters and given volume fractions.
Allaire~\etal~\cite{AlGePa19,GeAlPa20} investigated the shape and topology optimization of a perforated material, which is macroscopically modulated and oriented.
To this end, they first computed the homogenized properties of a certain class of microstructures. 
Then, they optimized the homogenized elastic parameters in the sense of a restricted free material optimization problem.
In the end, they properly placed the resulting local material pattern in a globally consistent way.
Schury~\etal~\cite{ScStWe12} used inverse homogenization to generate microstructures.
They solved the macroscopic problem with a free material optimization procedure. 
In a post-processing step they incorporated connectivity constraints.
Valentin~\etal~\cite{VaHuSt20} implemented hierachical B-splines on sparse grids to interpolate Cholesky factors of the parametrized materials on the microcells.
Groen~\etal~\cite{GrSi18} proposed a projection of rotationally symmetric spatially varying microstructures to generate two-dimensional lattice materials.
Ferrer~\etal~\cite{FeCaHe18} precomputed in 2D an offline dictionary of optimal microstructures \wrt input stresses represented on the two-dimensional sphere.
Zhou~\etal~\cite{ZhDuKi19} applied a shape metamorphosis to properly connect adjacent microcells.
Finally, we refer to the recent review article by Wu~\etal~\cite{WuSiGr21} about topology optimization of multi-scale structures.

\paragraph*{Hashin--Shtrikman bounds.} 
If we only consider homogenized elasticity tensors which are isotropic, 
the well-known Hashin--Shtrikman bounds~\cite{HaSh63} provide an estimate for the set of achievable pairs of parameters.
Precisely, all possible pairs of Poisson's ratio $\nu$ and Young's modulus $E$ are situated in a triangle with base $[0,1]$ on the $\nu$ axes. 
To achieve arbitrary pairs $(\nu,E)$ in this triangle, 
nested laminated microstructures have to be taken into account~\cite{MiBrHa17}.
As far as we know, the precise subset of material parameters which can be attained by a microstructure inside the triangle given by the Hashin--Shtrikman bounds is unknown.
Recently, in \cite{OsOvTo18}, specific microstructures have been studied to explore the Hashin--Shtrikman bounds.
More precisely, the authors verified numerically that a specific class of microstructures
exhausts a large range of the possible isotropic homogenized elasticity tensors. 
In particular, for zero stiffness of the soft material and positive Poisson's ratio, 
these microstructures experimentally seem to allow a fairly close approximation of the Hashin--Shtrikman bounds.
In \cite{MePoTo19}, certain microstructure symmetries are shown to imply symmetries of the homogenized elasticity tensor. 
Furthermore, crystal symmetries have been used to study Hashin--Shtrikman bounds in \cite{YeRoMe20}.

\paragraph*{Shape optimization tools.}
For shape and topology optimization, the implicit representation of shapes via level sets proved to be very powerful 
as investigated by Allaire and coworkers \cite{Al04,AlBoFr97,AlGoJo04,AlJoTo02}
as well as Wang~\etal~\cite{WaWaGu03}.
The homogenization method was used in the seminal paper by Bends{\o}e~\etal~\cite{BeKi88} 
for the numerical solution of an elastic shape optimization problem.
Optimal microstructures have been explored by Guedes~\etal~\cite{GuRoBe03}
to achieve approximative bounds on the elastic energy.
An alternating two-scale optimization was proposed by Barbarosie~\etal~\cite{BaTo14}. 
In \cite{AlDa14}, Allaire~\etal studied the optimization of multi-phase materials with a regularization of the shape derivative in the level set context. 
In this paper, we consider an implicit description of a microscopic shape via a phase-field approach using the double well approach by Modica and Mortola~\cite{MoMo77}. 
This approach has been employed in elastic shape optimization by Wang and Zhou~\cite{ZhWa07}, Blanck~\etal~\cite{BlGaHa14, BlGaHe16} using linearized elasticity, and in \cite{PeRuWi12} in the context of nonlinear elasticity. 
Guo~\etal~\cite{GuZhWa05} described the characteristic function of the object by the composition of a smoothed Heaviside function with a level-set function, where the smoothed Heaviside function acts like a phase-field profile. 
Wei and Wang~\cite{WeWa09} encoded the object via a piecewise constant level set function closely related to the phase-field approach.
An equality constraint for the volume can either be ensured by a Cahn--Hilliard-type $H^{-1}$-gradient flow \cite{ZhWa07}
or a Lagrange-multiplier ansatz \cite{AlJoTo03,BoCh03,LiKoHu05}. 
An inequality constraint for the volume of the elastic object has been implemented in \cite{WaWaGu03, WaZhDi04, GuZhWa05} using a Lagrange multiplier approach.
In general, shape optimization problems of the above type are ill-posed, since minimizing sequences of characteristic functions associated with elastic objects do not necessarily converge in a weak topology to an elastic object. 
We refer to \cite[Theorem 4.1.2]{Al02} for a detailed treatment of this point.
To avoid this ill-posedness, the void phase is typically replaced by some weak material and a perimeter penalty is added to the objective functional (\cf~\cite{AmBu93} for a scalar problem).  
In \cite{StKoLe09}, Stingl \etal proposed a convex semidefinite programming method for a free material optimization,
where inequality constraints on the elasticity tensor are taken into account.
Finally, we refer to the recent survey article \cite{AlDaJo21}.

\bigskip

\paragraph*{Organization.}
The paper is organized as follows.
In Section~\ref{sec:twoscale}, we review the modeling of two-scale elasticity.
This in particular includes the formulation of a corrector problem including a manufacturability constraint.
This corrector problem has to be solved to identify the homogenized macroscopic elasticity tensor for a given periodic microscopic material pattern.
Section~\ref{sec:micro-scale} investigates the microscopic shape and topology optimization using a phase-field approximation.
We show numerical results of the offline fine-scale optimization in Section~\ref{sec:fineScaleResults}.
Then, in Section~\ref{sec:macroopt} the 
restricted 
free material optimization approach on the macro-scale is discussed.
Here, we describe how to define a regular parameterization of admissible effective isotropic material parameters and associated cost values based on a strategy for the  sampling of elasticity tensors.
Furthermore, we discuss how to realize an approximation of a workpiece with optimal homogenized macroscopic elastic material properties with an additively manufactured fine-scale pattern based on a specific isotropic elastic material.
Finally, we present numerical results of macroscopically optimized workpieces with additively manufacturable fine-scale structures in Section~\ref{sec:twoScaleResults}.

\section{A two-scale elasticity model with manufacturability constraint} \label{sec:twoscale}
Let us consider  a bounded macroscopic domain 
$\domain \subset \R^d$ ($d=2,3$)
describing an elastic workpiece 
with total free energy $\Etot[\Etensor^\epsilon,\displaceMacro] \coloneqq \Eelast[\Etensor^\epsilon,\displaceMacro] - \Epot[\displaceMacro]$, where
\begin{align}\label{eq:elastenergy}
\Eelast[\Etensor^\epsilon,\displaceMacro] \coloneqq \tfrac12 \int_\domain \Etensor^\epsilon \varepsilon[\displaceMacro] : \varepsilon[\displaceMacro] \d x \,, \quad 
\Epot[\displaceMacro] \coloneqq \int_\domain f \cdot \displaceMacro \d x + \int_{\Gamma_N}  g \cdot \displaceMacro \d a
\end{align}
are the stored elastic energy of linearized elasticity and the potential energy, respectively, and 
$\varepsilon(\displaceMacro) \coloneqq \tfrac12 (D\displaceMacro^T+D\displaceMacro)$ is the strain tensor 
with $D\displaceMacro$ denoting the Jacobian (the weak differential) of the displacement $\displaceMacro \in \displaceMacro^\partial + H^{1,2}_{\Gamma_D} (\domain)^d$ and 
$\Etensor \,A:B \coloneqq \sum_{i,j,k,l} \Etensor_{ijkl} A_{ij} B_{kl}$ for $A,B \in \R^{d\times d}$.
The boundary $\partial \domain$ is split into a Dirichlet boundary $\Gamma_D$ and a Neumann boundary $\Gamma_N$, 
$H^{1,2}_{\Gamma_D} (\domain)^d$ is the space of displacements with square integrable Jacobian which vanish on $\Gamma_D$, and $U^\partial$ is a fixed extension of sufficiently smooth Dirichlet data.

The elasticity tensor is assumed to depend on the geometry at a fine-scale $\epsilon> 0$,
\begin{align}
\Etensor^\epsilon(x) =  \Etensor\left(x,\tfrac{x}{\epsilon}\right)
\end{align}
with $\Etensor(x,y) = \chi(x,y) \Etensor^1 + (1-\chi(x,y)) \Etensor^2$,
where $\chi \colon \domain \times \R^d \to \{0,1\};\ (x,y) \mapsto \chi(x,y)$ with $\chi$ being $Q$-periodic in $y$ with $Q\coloneqq[0,1]^d$,
which reads as $\chi(x,y+\Z^d) = \chi(x,y)$.
In fact, $\chi(x,\cdot)$ is considered as a characteristic function describing a local fine-scale geometry 
on the unit cube $Q$ and the two tensors $\Etensor^1$ and $\Etensor^2$ are the elasticity tensors on the two material phases
identified by $\{\chi(x,\cdot)=1\}$ and $\{\chi(x,\cdot)=0\}$, respectively.
In our application, we assume for simplicity that the elasticity on the fine-scale 
is of Lam{\'e}-Navier type with a hard and a soft phase, \ie
\begin{align}\label{eq:isotropC}
 \Etensor^1  \coloneqq 
 \left( \kappa \delta_{ij} \delta_{kl} 
      + \mu \left(\delta_{ik}  \delta_{jl} +\delta_{il} \delta_{kj} - \tfrac{2}{d} \delta_{ij} \delta_{kl} \right) \right)_{i,j,k,l=1,\ldots, d} 
      \,, \quad 
 \Etensor^2 \coloneqq \delta \Etensor^1
\end{align}
where $\mu$ denotes the shear modulus, $\kappa$ the bulk modulus and 
$0 < \delta \ll 1$.
For our application it will be more convenient to take into account Poisson's ratio $\nu$
and Young's modulus $E$,
which are given by
\begin{align}
  \nu = \tfrac{\kappa - \mu}{\kappa + \mu}, \quad E = \tfrac{4 \kappa \mu}{\kappa + \mu} \quad \text{for } d=2 \ \text{ and }   
 \nu=\tfrac{3\kappa - 2\mu}{6 \kappa + 2\mu}, \quad E=\tfrac{9\kappa \mu}{3 \kappa + \mu} \quad \text{for } d=3 \, .
\end{align}
In the entire paper we assume that $\Etensor^1$ and $\Etensor^2$ are both bounded and coercive 
(which ensures well-posedness of the elastic problem since $\delta > 0$).
Note that the case $\delta = 0$ would correspond to a pure void material, 
which is replaced by a very soft material with relative stiffness $\delta \ll 1$.

\paragraph*{Effective, homogenized macroscopic elasticity.}
It is known from homogenization theory that in the limit $\epsilon \to 0$,
a linearized elasticity model is observed
with an effective elasticity tensor field $\Etensor^\ast \colon \domain \to \R^{d^4}$
with $\Etensor^\ast$ depending on the material microstructure described by $\chi(x, \cdot)$.
We write $\Etensor^\ast(x) = \Etensor^\ast[\chi(x,\cdot)]$ for any $\chi(x,\cdot) \in L^\infty(Q;\{0,1\})$.
The effective elasticity tensor field $\Etensor^\ast$ is bounded and coercive,
fulfills the symmetries $\Etensor^\ast_{ijkl} = \Etensor^\ast_{jikl}  = \Etensor^\ast_{ijlk} = \Etensor^\ast_{klij}$ 
for all $i,\ j,\ k,\ l = 1,\ldots, d$, and is characterized by
\begin{align}
 C^\ast_{ijkl}(x) 
 &= \int_Q C(x,y) \left(\delta_{kl} + \varepsilon[\displaceMicro^{kl}(y)]\right) : \left(\delta_{ij} + \varepsilon[\displaceMicro^{ij}(y)]\right) \d y
\end{align}
with $\displaceMicro^{kl}\in H^{1,2}_{\sharp}(Q)^d$ solving the corrector problem given by
\begin{align}
 0 =  \int_Q \sum_{i,j,m,n=1,\ldots, d} C_{ijmn}(x,y) \left(\delta_{mk} \delta_{nl} + \varepsilon_{mn}[\displaceMicro^{kl}(y)] \right) \varepsilon_{ij}[\phi (y)] \d y 
\end{align}
for all $\phi \in H^{1,2}_{\sharp}(Q)^d$.
Here, $H^{1,2}_{\sharp}(Q)^d =\{ \displaceMicro \in H^{1,2}_{loc}(\R^d)^d\ \vert \ \int_Q \displaceMicro(y) \d y =0\,,\ \displaceMicro(y+\mathbb{Z}^d) = \displaceMicro(y) \}$
is the set of displacements in the Sobolev space of locally square integrable Jacobian, which are $Q$ periodic and have mean value zero.
Further, $\Etensor^\ast$ depends continuously on $\chi$, with respect to the strong $L^1(Q)$ norm \cite[Lemma~2.1.6]{Al02}.

We remark that the case of a hard and a void phase ($\delta = 0$) is more subtle.
Indeed, the microscopic elastic domains $\{y\in Q\ \vert \ \chi(y) = 1\}$ might be disconnected and thus the corrector problem may not be solvable.
We shall not address this case in the present paper.

\begin{figure}[htbp]
\centering{ 
 \includegraphics[width=0.25\textwidth]{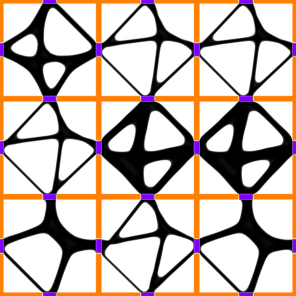} 
 \hspace{5ex} 
 \includegraphics[width=0.25\textwidth]{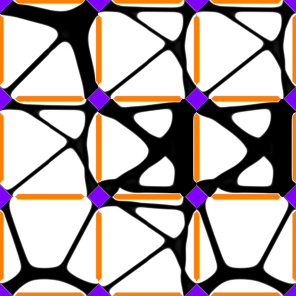}
 \hspace{5ex} 
 \includegraphics[width=0.25\textwidth]{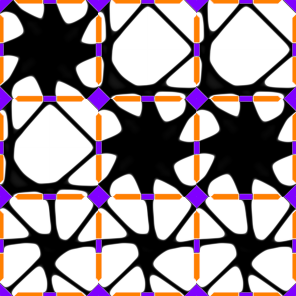}
}
 \caption{
   A sketch of three different manufacturability constraints on $Q=[0,1]^2$ 
   with the sets $B^1$ (violet) and $B^2$ (orange).
   For each of them a three-by-three pattern of fine scale cells is displayed.
 }
 \label{fig:bridges}
\end{figure}

\paragraph*{Manufacturability constraint.} 
We take into account the additive manufacturability by enforcing material 
fixed
bridges at the boundary of the microscopic cubic domain $Q$.
Assuming symmetries of these bridges at opposite faces, 
this guarantees that different microstructures fit together at the boundaries.
More precisely, we consider so-called bridge sets $B^1$ and $B^2$, which are assumed to be open 
Lipschitz subsets of $Q$ with
\begin{align}
 \partial Q \subset \partial {B^1} \cup \partial {B^2},\ B^1 \cap B^2 = \emptyset,\  
 \mathcal{H}^{d-1}(\partial B^i\cap \partial Q ) > 0 \text{ for $i=1,2$.}
\end{align}
We stress that these sets are fixed and do not depend on the choice of the microstructure in the individual cell.
Furthermore, we assume that a bridge domain attached to one face of the cubic cell has a counterpart on the opposite face of the cell,
which coincides due to the periodicity of $\chi$ to the bridge structure in the neighbouring cell adjacent to the face, 
\ie we require
\begin{align}
y \in \partial B^i \cap \{y_k=1\} \Leftrightarrow  y - e_k \in \partial B^i 
\end{align}
with $e_k$ being the $k$th canonical basis vector in $\R^d$.
Practically relevant choices of the sets $B^i$ have a simple form.
In Figure~\ref{fig:bridges} we show three different examples for the bridge sets $B^1$ and $B^2$, which we will study later numerically.

Finally, we are in the position to define the set of characteristic functions representing a manufacturable microstructure
\wrt the fixed bridge sets $B^1$ and $B^2$ by
\begin{align}
 \mathcal{B} \coloneqq \{  \chi \in BV(Q,\{0,1\})  \ \vert \   \chi(x) = 1 \ \text{ for a.e. } x \in B^1, \chi(x) = 0 \ \text{ for a.e. } x \in B^2  \} \, .
\end{align}

We remark that other constraints are actually relevant to guarantee manufacturability,
which are not explicitly enforced in our model.
Due to the perimeter regularization, the feature size is implicitly not allowed to become arbitrarily small.
Overhangs, which we will indeed observe in our numerical results, are not forbidden.

\section{Shape optimization on the micro-scale} \label{sec:micro-scale}
In Section~\ref{sec:twoscale} we defined the effective elasticity tensors $\Etensor^\ast$ 
reflecting the material pattern described by a characteristic function on the fundamental domain $Q$ on the micro-scale.
Now, we ask for which elasticity tensors $\Etensor$
there exists a microscopic material pattern in the class of manufacturable 
patterns with this tensor as the effective elasticity tensor,
in the sense that $\Etensor=\Etensor^\ast[\chi]$. 
Therefore, we aim to determine the underlying microscopic shape $\chi$ as the solution of an optimization problem for a specific cost functional
of an elastic shape optimization problem on the periodic unit cell.
We will first run this shape optimization problem in an offline phase for a sufficiently large set of elasticity tensors, and then retrieve via interpolation suitable approximations for the set of realizable elasticity tensors and for the value of the cost functional on this set.
This will then be used in an online phase for the free material optimization given a particular load scenario on the macro-scale.

To describe the material pattern on the macro-scale, 
we first define for a given effective elasticity tensor $\Etensor$ 
the set of characteristic functions $\chi$
representing a microscopic material pattern corresponding to $\Etensor$ by
\begin{align}
 \left\{ \chi  \in \mathcal{B} \  \Big\vert \  \Etensor=\Etensor^\ast[\chi] \right\}\,.
\end{align}
Note that this set might be empty, which means that the tensor $\Etensor$ cannot be recovered by a characteristic function in $\mathcal{B}$.
However, continuity of $\Etensor^\ast$ and the structure of $\mathcal B$ show that this set is closed with respect to the strong $L^1(Q)$ topology, and hence also with respect to the weak $BV(Q)$ topology.

Now, on the micro-scale, we optimize the material distribution over this set.
To this end, we introduce the associated microscopic cost functional  
\begin{align} 
 \mathbf{J}_{\text{micro}}[\chi]  \coloneqq  c_V \Volume[\chi] + c_P \Perimeter[\chi] 
 = c_V \int_Q \chi(y) \d y + c_P |D\chi |(Q) \,,
\end{align}
where $\Volume[\chi] \coloneqq \int_Q \chi \d y$ is the volume of $\{ y \in Q \ \vert \ \chi(y) = 1 \}$,
$\Perimeter[\chi] \coloneqq |D\chi|(Q)$ is the perimeter given by the total variation of $\chi$ on $Q$, 
and $c_V \geq 0$, $c_P > 0$ are corresponding microscopic weights.
Let us remark that $c_P > 0$ is needed to ensure the well-posedness of the variational problem 
(\cf~\cite{PeRuWi12} for a proof in the case of nonlinear elasticity).
In the case $c_P=0$, we do not expect existence of minimizers, as minimizing sequences only converge weakly in $L^p(Q)$, and the constraint $\chi(y)\in\{0,1\}$ a.e. is lost in the limit.
Next, we select one (arbitrary) minimizer of $\mathbf{J}_{\text{micro}}$ 
\begin{align}\label{eq:shapeOptProblemCell}
 \chi^\ast[\Etensor] \in 
 \argmin_{\chi \in \mathcal{B} \,,\ \Etensor=\Etensor^\ast[\chi]} \mathbf{J}_{\text{micro}}[\chi] \,.
\end{align}
Finally, we introduce a cost function $\mathbf{j}[\cdot]$ which gives the value of $\mathbf{J}_{\text{micro}}$ on the microstructure $\chi$ chosen for a given $\Etensor$.
It takes into account the local mass and perimeter density which are required to realize a specific elasticity tensor on the macroscopic domain. 
This value will be used as a cost density in the free material optimization approach on the macro-scale in Section~\ref{sec:macroopt}. 
Precisely,
\begin{align}\label{eq:localCost}
 \mathbf{j}[\Etensor] 
 \coloneqq \mathbf{J}_{\text{micro}}[\chi^\ast[\Etensor]] 
 = c_V \int_Q \chi^\ast[\Etensor](y) \d y + c_P |D \chi^\ast[\Etensor] |(Q)
\end{align}
for $\Etensor \in \mathcal{C}_{\text{real}}$,
and $\mathbf{j}[\Etensor] = \infty$ for $\Etensor \notin \mathcal{C}_{\text{real}}$,
where $\mathcal{C}_{\text{real}}$ denotes the set of realizable effective elasticity tensors.
Here, in the evaluation of the cost functional $\mathbf{J}_{\text{micro}}$ and for the purpose of the definition of  $\mathbf{j}[\Etensor]$ 
we replace the weight $c_P>0$  by a possibly different weight $\widehat{c_P} \ge 0$,
because of the different scaling of the volume with $\epsilon^d$ and the perimeter with $\epsilon^{d-1}$.
This does not alter the definition of $\chi^\ast[\Etensor]$ in \eqref{eq:shapeOptProblemCell}, therefore
$\mathbf{j}[\Etensor]$ reflects a possibly modified cost value for the minimum of $\mathbf{J}_{\text{micro}}[\chi]$.
Note that the method proposed here could be applied to recover any tensor $\Etensor^\ast$ 
on the full space of all effective elasticity tensors,
which is six-dimensional in 2d and 21-dimensional in 3d.
For computational simplicity, in our application, we will consider $\mathcal{C}_{\text{adm}}$ as the (two dimensional) space of isotropic elasticity tensors, 
which can be parametrized over the bulk and shear modulus.

\paragraph*{Phase-field approximation.}
A direct numerical treatment of the characteristic function $\chi \in BV(Q, \{0,1\})$ is algorithmically quite demanding,
therefore we replace it by a phase-field function 
$\pf \in H^{1,2}(Q,\R)$
allowing intermediate values in $(-1,1)$ that are appropriately penalized.
Correspondingly, the perimeter $|D\chi|(Q)$ is replaced by the Modica--Mortola type functional 
\begin{align}
 L^{\sigma}[\pf] 
 \coloneqq \tfrac12 \int_Q \sigma|\nabla \pf|^2 + \tfrac1\sigma \DoubleWell(\pf) \d y \,,
\end{align}
where $\sigma > 0$ describes the width of the diffuse interface between the soft and hard material phases (\cf~\cite{PeRuWi12})
and $\DoubleWell(\pf) \coloneqq \tfrac{9}{16}(\pf^2-1)^2$ has
two minima at $\pf=-1$ and $\pf=1$.
In the limit $\sigma\to 0$, the phase-field $\pf$ is expected to lead to a clear separation between two pure phases $-1$ and $1$  
and $L^{\sigma}$ is known to $\Gamma$-converge to the perimeter functional, defined as  $\tfrac12|Dv|(Q)$ for $v\in BV(Q;\{-1,1\})$, and $\infty$ otherwise, see \cite{Br02}.
We correspondingly adapt the notion of the phase dependent elasticity tensor
$\Etensor[\pf](y) = \chi[\pf(y)] \Etensor^1 + (1-\chi[\pf(y)]) \Etensor^2$.
Here the characteristic function $\chi$ has been replaced by its approximation via the phase-field variable $v$, which we take to be 
$\chi[\pf] \coloneqq \tfrac{1}{16} (1 + \pf)^4$.
Note that due to the $\Gamma$-convergence of $L^{\sigma}$, the exact choice of $\chi[\pf]$ does not matter in the limit $\sigma \to 0$ as long as $\chi[-1]=0$ and $\chi[1]=1$.
Then the stored elastic energy is redefined in terms of the phase-field function $\pf$ and 
the displacement $\displaceMicro$ as
\begin{align}
 \Eelast[\pf,\displaceMicro] 
 \coloneqq 
 \tfrac12 \int_Q  \Etensor[\pf] \varepsilon(\displaceMicro) : \varepsilon(\displaceMicro) \d y.
\end{align}
Now, we obtain the phase-field adapted effective elasticity tensor 
\begin{align}
C^\ast_{ijkl}[v] 
 &= \int_Q C[\pf](y) \left(\delta_{kl} + \varepsilon[\displaceMicro^{kl}(y)]\right) : \left(\delta_{ij} + \varepsilon[\displaceMicro^{ij}(y)]\right) \d y
\end{align}
where $\displaceMicro^{kl}\in H^{1,2}_{\sharp}(Q)^d$ solves the corrector problem
\begin{align}
 0 = \int_Q  \sum_{i,j,m,n=1,\ldots, d} C_{ijmn}[\pf](y)\left(\delta_{mk}  \delta_{nl} + \varepsilon_{mn}[\displaceMicro^{kl}(y)]\right) \varepsilon_{ij}[\phi(y)] \d y 
\end{align}
for all $\phi \in H^{1,2}_{\sharp}(Q)^d$.
Next, to encode the printability constraint given by the sets $B^1$ and $B^2$, we remark that a phase-field function $\pf \in H^{1,2}(Q)$ is not allowed to jump.
Instead, $\pf$ should perform a transition between $B^1$ and $B^2$ on a length scale $\sigma$.
Thus, we define the set
\begin{align}
\mathcal{B}_{\sigma} \coloneqq \{  \pf \in H^{1,2}(Q)  \ \vert \  \pf = 1 \text{ on } B^1, \; \pf = -1 \text{ on } B^2 \setminus B_{\sigma}(B^1) \} \, ,
\end{align}
where $B_{\sigma}(B^1) \coloneqq \bigcup_{x \in B^1} B_{\sigma}(x)$ denotes the 
set of points whose distance from the set $B^1$ is smaller than $\sigma$.
Finally, we define the optimal phase-field $\pf^\ast[\Etensor]$ realizing an effective elasticity tensor $\Etensor$ as 
\begin{align} 
 \pf^{\ast, \sigma}[\Etensor]   \coloneqq
 \argmin_{\pf \in   \mathcal{B}_{\sigma} \,, \ \Etensor=\Etensor^\ast[\pf]}\left( c_V \int_Q \tfrac{1}{16} \left(1 + \pf(y)\right)^4 \d y + c_P L^{\sigma}[\pf]\right) \,.
 \end{align}
The microscopic cost density used for the macroscopic optimization is defined as 
\begin{align}
\mathbf{j}^\sigma[\Etensor] \coloneqq c_V \int_Q \tfrac{1}{16} \left(1 + \pf^{\ast, \sigma}[\Etensor](y)\right)^4 \d y + 
\widehat{ c_P }
L^{\sigma}[\pf^{\ast,\sigma}[\Etensor]] \,,
\end{align}
where $\widehat{c_P} \geq 0$.
Let us remark that in case $\widehat{c_P}=0$ the macroscopic optimization does not directly depend on the 
microscopic perimeter but only indirectly via the choice of $\pf^{\ast,\sigma}$.

\textit{Remark (on the $\Gamma$-convergence of $\mathbf{j}^\sigma$):}
As discussed above, the functionals $L^{\sigma}$ $\Gamma-$converge to the area functional with respect to the strong $L^1$ convergence \cite{MoMo77,Br02}. The growth of $W$ implies that 
$\int_{|v_\sigma|\ge 2} |v_\sigma|^4 dy\to0$ 
along sequences with bounded energy,
so that
the first term is continuous, and 
one can easily deduce that the functionals
\begin{align}\label{eqgammalsig}
 c_V \int_Q \tfrac{1}{16} \left(1 + \pf(y)\right)^4 \d y + c_P L^{\sigma}[\pf]
\end{align}
$\Gamma-$converge, in the limit $\sigma\to0$, to the sharp-interface limit defined by
\begin{align}\label{eqgammal0}
 c_V \int_Q \tfrac{1}{16} \left(1 + \pf(y)\right)^4 \d y + \frac12c_P |D\pf|(Q)
\end{align}
if $\pf\in BV(Q;\{\pm1\})$, and $\infty$ otherwise.
Here, we are interested in the restrictions of these functionals to the space 
$X_\Etensor:=\{ \pf : \ \Etensor=\Etensor^\ast[\pf]\}$
by setting them equal to $\infty$ on the complement.
The space $X_\Etensor$
is closed with respect to the strong $L^1$ convergence
(see \cite[Prop.~1.4.5]{Al02}), therefore the $\Gamma-\liminf$ inequality survives the restriction without any change. The situation is, however, more complex for the $\Gamma-\limsup$ inequality, as such restrictions pose constraints that the recovery sequence must obey. Indeed, already in the  (simpler) case of the gradient theory of phase transitions, if one restricts $L^\sigma$ to those test function $\pf$ that obey a certain volume constraint (for example prescribing the value of  $\int \pf \d y$), then special care  needs to be taken both in the density result and in the explicit construction of a recovery sequence in order to avoid violating the constraint \cite[Lemma~1 and Prop.~2]{Mo87}. The present case in which one prescribes the value of all components of $\Etensor^\ast[\pf]$ is more complex: firstly, one does not prescribe a single scalar quantity but several, and secondly, these quantities do not depend linearly on $\pf$.
We assume here that for the relevant values of $\Etensor$
the space $X_\Etensor$ is either empty or sufficiently rich that a recovery sequence can be constructed. 
We expect this to be true if $\Etensor$ is not on the boundary of the admissible set. This assumption in turn implies that the restriction of \eqref{eqgammalsig} to $X_\Etensor$
$\Gamma$-converges to the restriction of \eqref{eqgammal0} to the same space. Similar comments apply (with heavier notation, but essentially the same difficulty) if one additionally restricts to $v\in \mathcal{B}_{\sigma}$.
Under these assumptions, for $\sigma\to0$ we have
(up to subsequences, and up to the choice of minimizer in \eqref{eq:shapeOptProblemCell})
$\chi[v^{\ast, \sigma}]\to \chi^\ast[\Etensor]$
in $L^p(Q)$, for any $p\in[1,\infty)$,
and finally
(as long as $\widehat{c_P}=c_P$)
$\mathbf{j}^\sigma[\Etensor]\to\mathbf j[\Etensor]$.
In our computations we use $\mathbf{j}^\sigma[\Etensor]$ as a cost density on the macro-scale.

\paragraph*{Finite element discretization.}
To discretize our computational domain, which is always given by the unit cube $Q = [0,1]^d$,
we use in the case $d=2$ a uniform quadrilateral mesh
and in the case $d=3$ a uniform hexahedral mesh.
In both cases, we denote by $h$ is the corresponding mesh size,
$\mathcal{K}_h$ the set of (quadrilateral or hexahedral) cells,
and $\mathcal{N}_h$ the set of nodes of the mesh. 
Then we define the space $\mathcal{V}_h$ of piecewise $d$-affine, continuous functions.
We consider discrete phase-fields $\pf_h \in \mathcal{V}_h$. 
For an index $\alpha \in \mathcal{N}_h$, we denote by $\phi_h^\alpha \in \mathcal{V}_{h}$ the nodal valued base function with $\phi_h^\alpha(\gamma) = \delta_{\alpha}(\gamma)$ for all $\gamma \in \mathcal{N}_h$.
Furthermore, we denote by 
$\overrightarrow{\pf}_h 
= \left( (\overrightarrow{\pf}_h)^\alpha \right)_{\alpha \in \mathcal{N}_h}
= \left( \pf_h(\alpha) \right)_{\alpha \in \mathcal{N}_h}
$ 
the representation of $\pf_h$ at nodal values.
Then the admissible set of discrete phase-fields taking into account the material bridges is given by
\begin{align}
 \mathcal{B}_{\sigma,h} \coloneqq \{  
    \pf_h \in \mathcal{V}_{h}  \ \vert \  
      \pf_h(\alpha) = 1 \text{ for all } \alpha \in B^1, \; 
      \pf_h(\alpha) = -1 \text{ for all } \alpha \in B^2 \setminus B_{\sigma}(B^1) 
    \} \, .
\end{align}
Moreover, we consider discrete correctors $\displaceMicro^{kl}_h \in \mathcal{V}_{h}^d$ and use the same notation $\overrightarrow{\displaceMicro_h}^{kl}_h$ for the representation at nodal values.
To implement fully discrete functionals, we apply a tensor product Simpson quadrature rule.
For each cell $K \in \mathcal{K}_h$, we denote by $\mathcal{S}$ these quadrature points, $\omega(s)$ the quadrature weight for a single quadrature point $s \in \mathcal{S}$, and $y(K,s)$ the corresponding mapping to coordinates in $Q$.
This allows to define a discrete version of the stored elastic energy by 
\begin{align}
  E_{\text{elast},h}[\pf_h,\displaceMicro_h] 
  & \coloneqq 
  \tfrac12 \sum_{K \in \mathcal{K}_h} |K| \sum_{s \in \mathcal{S}} \omega(s) \,  \Etensor[\pf_h(y(K,s))] \varepsilon(\displaceMicro_h(y(K,s))) : \varepsilon(\displaceMicro_h(y(K,s))).
\end{align}
First, for a fixed phase-field $\pf_h$, we describe the numerical scheme to compute a solution of the corrector problem.
To this end, we define a system matrix $L^{kl}_h[\pf_h] = \left(L^{kl}_{\alpha,\beta}[\pf_h]\right)_{\alpha,\beta \in \mathcal{N}_h^d}$
and a right hand side $R^{kl}_h[\pf_h] = \left(R^{kl}_\alpha[\pf_h]\right)_{\alpha \in \mathcal{N}_h^d}$ by
\begin{align*}
 L^{kl}_{\alpha,\beta}[\pf_h] & \coloneqq 
 \sum_{K \in \mathcal{K}_h} |K| \sum_{s \in \mathcal{S}} \omega(s) \,  
 \sum_{i,j,m,n=1,\ldots, d} C_{ijmn}[\pf_h](y(K,s)) \varepsilon_{mn}[\phi_h^{\beta}(y(K,s))] \varepsilon_{ij}[\phi_h^{\alpha}(y(K,s))] \, , \\
 R^{kl}_\alpha[\pf_h] & \coloneqq - \sum_{K \in \mathcal{K}_h} |K| \sum_{s \in \mathcal{S}} \omega(s) \,  
 \sum_{i,j,m,n=1,\ldots, d} C_{ijmn}[\pf_h](y(K,s)) \delta_{mk} \delta_{nl} \varepsilon_{ij}[\phi_h^{\alpha}(y(K,s))] \, .
\end{align*}
After applying the periodic boundary condition to $L^{kl}_h$ and $R^{kl}_h$ and incorporating the mean value constraint $\int_Q \displaceMicro^{kl}_h \d y = 0$,
the solution to the corrector problem $\displaceMicro^{kl}_h[\pf_h]$ explicitly depending on $\pf_h$
is obtained by the linear system $L^{kl}_h[\pf_h] \overrightarrow{\displaceMicro_h}^{kl}_h[\pf_h] = R^{kl}_h[\pf_h]$.
This allows to formulate a fully discrete version of the constraint $\Etensor = \Etensor^{\ast}[\chi]$
by functions 
\begin{align}
 g_{ijkl}^{\Etensor}[\pf_h] = \hat g_{ijkl}^{\Etensor}[\pf_h, \displaceMicro^{kl}_h[\pf_h], \displaceMicro^{ij}_h[\pf_h] ] \, ,
\end{align}
where $\hat g_{ijkl}^{\Etensor}$ explicitly depending on 
$\pf_h$, $\displaceMicro^{kl}_h$, and $\displaceMicro^{ij}_h$ is defined by
\begin{align}
 \hat g_{ijkl}^{\Etensor}[\pf_h, \displaceMicro^{kl}_h, \displaceMicro^{ij}_h] \coloneqq
 \sum_{K \in \mathcal{K}_h} |K| \sum_{s \in \mathcal{S}} \omega(s) \, 
 C[\pf_h](y(K,s))
 \big( \delta_{kl} + \varepsilon[\displaceMicro^{kl}_h(y(K,s))] \big) :
 \big( \delta_{ij} + \varepsilon[\displaceMicro^{ij}_h(y(K,s))] \big)
 - \Etensor_{ijkl}  \, .
\end{align}
Thus, we require $g_{ijkl}^{\Etensor}[\pf_h] = 0$ for all $i,j,k,l$.
By symmetries of the effective elasticity tensor, it is sufficient to enforce these constraints only on the corresponding Voigt tensor,
\ie we end up with six constraints in 2d and $21$ constraints in 3d.

\medskip

Now, a discrete versions of the microscopic cost functional is given by
\begin{align}
  \mathbf{J}_{\text{micro},h}[\pf_h] 
  & \coloneqq \sum_{K \in \mathcal{K}_h} |K| \sum_{s \in \mathcal{S}} \omega(s) \, 
  \left( \tfrac{c_V }{16} (1+\pf_h(y(K,s)))^4 
   + \tfrac{c_P}{2} \left( \sigma |\nabla \pf_h(y(K,s))|^2 + \tfrac1\sigma \DoubleWell(\pf_h(y(K,s))) \right)
   \right) \,.
\end{align}
Finally, we can formulate a fully discrete version of the shape optimization problem~\eqref{eq:shapeOptProblemCell} by
\begin{align} \label{eq:shapeOptProblemCellFEM}
  \pf_h^\ast[\Etensor] \in 
 \argmin_{\pf_h \in \mathcal{B}_{\sigma,h} \,,\ g_{ijkl}^{\Etensor}[\pf_h] = 0} \mathbf{J}_{\text{micro},h}[\pf_h] \,.
\end{align}
For the numerical solution, we use the quasi-Newton BFGS method from the IPOPT package \cite{WaBi06}, 
which provides a general interface to solve constraint minimization problems of the above type.
Here, to solve \eqref{eq:shapeOptProblemCellFEM}, we also have to provide the derivatives of 
$\mathbf{J}_{\text{micro},h}$ and $g_{ijkl}^{\Etensor}$.
Moreover, note that the condition $\pf_h \in \mathcal{B}_{\sigma,h}$ turns into simple nodal-valued box constraints.
Here, the derivative of $g_{ijkl}[\pf_h]$ is given by
\begin{align}
 \frac{d}{d \pf_h} g_{ijkl}^{\Etensor}[\pf_h] 
&= \partial_{\pf_h}  \hat g_{ijkl}^{\Etensor}[\pf_h, \displaceMicro^{kl}_h[\pf_h],\displaceMicro^{ij}_h[\pf_h]] \\
& + \partial_{\displaceMicro^{kl}_h}  \hat g_{ijkl}^{\Etensor}[\pf_h, \displaceMicro^{kl}_h[\pf_h],\displaceMicro^{ij}_h[\pf_h]] \, \partial_{\pf_h} \displaceMicro^{kl}_h[\pf_h] 
+ \partial_{\displaceMicro^{ij}_h}  \hat g_{ijkl}^{\Etensor}[\pf_h, \displaceMicro^{kl}_h[\pf_h],\displaceMicro^{ij}_h[\pf_h]] \, \partial_{\pf_h} \displaceMicro^{ij}_h[\pf_h] 
\, .
\end{align}
By the definition of the corrector problems, we observe for the partial derivatives that 
$\partial_{\displaceMicro^{kl}_h} \hat g_{ijkl}^{\Etensor}[\pf_h, \displaceMicro^{kl}_h[\pf_h],\displaceMicro^{ij}_h[\pf_h]] = 0$
and 
$\partial_{\displaceMicro^{ij}_h} \hat g_{ijkl}^{\Etensor}[\pf_h, \displaceMicro^{kl}_h[\pf_h],\displaceMicro^{ij}_h[\pf_h]] = 0$.
Thus, we finally obtain that 
\begin{align}
 \frac{d}{d \pf_h} g_{ijkl}^{\Etensor}[\pf_h] 
= \partial_{\pf_h}  \hat g_{ijkl}^{\Etensor}[\pf_h, \displaceMicro^{kl}_h[\pf_h], \displaceMicro^{ij}_h[\pf_h]] \, .
\end{align}

\section{Numerical results on the micro-scale} \label{sec:fineScaleResults}

In this section, we will present several numerical results for the shape optimization problem~\eqref{eq:shapeOptProblemCell}.
We will first discuss the case $d=2$.
For the discretization of the unit square $[0,1]^2 \subset \R^2$, we use a mesh size of $h=2^{-7}$
and $\sigma=2h$ for the interface width of the Modica--Mortola functional.
We always choose material parameters $\mu = 4$, $\kappa=6.67$ for the stiff phase,
which corresponds to a Poisson's ratio $\nu = 0.25$ and a Young's modulus $E = 10$.
Furthermore, we use $\delta = 10^{-4}$ as scaling factor for $\mu$ and $\kappa$ to approximate the void material in the weak phase.
Then, for a given volume fraction $\theta_s \in [0,1]$ of the stiff phase,
the upper Hashin--Shtrikman bounds 
(\cf \cite{HaSh63})
for the bulk and shear modulus   
\begin{align}
 \kappa_u = \kappa + \frac{1-\theta_s}{\frac{1}{(\delta-1)\kappa} + \frac{\theta_s}{\lambda + 2 \mu}} \,,
 \qquad
 \mu_u  = \mu + \frac{1-\theta_s}{\frac{1}{(\delta-1)\mu} + \frac{2 \theta_s (d-1) (\kappa + 2 \mu)}{ (d^2 + d - 2 ) \mu ( \lambda + 2 \mu ) }}
\end{align}
a priori allow to restrict the range of possible realizable isotropic effective elasticity tensors
to be located inside the triangle spanned by the points
$(\nu_{\text{min}},0)$,
$(\nu_{\text{max}},0)$,
and 
$(\nu_{\text{top}},E_{\text{top}})$
with 
\begin{align}
\nu_{\text{min}} = -1,
\qquad 
\nu_{\text{max}} =
\begin{cases}
 1 &\text{if } d = 2, \\
 \frac{1}{2} &\text{if } d = 3,
\end{cases}
\qquad
(\nu_{\text{top}},E_{\text{top}}) =
\begin{cases}
 \left( \frac{\kappa_u - \mu_u}{\kappa_u + \mu_u}, \; \frac{4\kappa_u \mu_u}{\kappa_u+\mu_u} \right)   & \text{if } d=2, \\
 \left( \frac{3\kappa_u - 2\mu_u}{2\kappa_u + 6\mu_u}, \; \frac{9\kappa_u \mu_u}{3\kappa_u+\mu_u} \right)  & \text{if } d=3.
\end{cases}
\end{align}
We shall consider the three  examples in Figure~\ref{fig:bridges} 
for the sets $B^1$ and $B^2$ that enforce the manufacturability constraint.
As stopping criterion for the IPOPT solver, we require that the overall residual is smaller than $10^{-10}$. 
The typical CPU time for one microstructure optimization is $10^3$ sec,
where points near the boundary usually converge slower than in the interior.

\begin{figure}[htbp]
\centering{
 \includegraphics[width=0.245\textwidth]{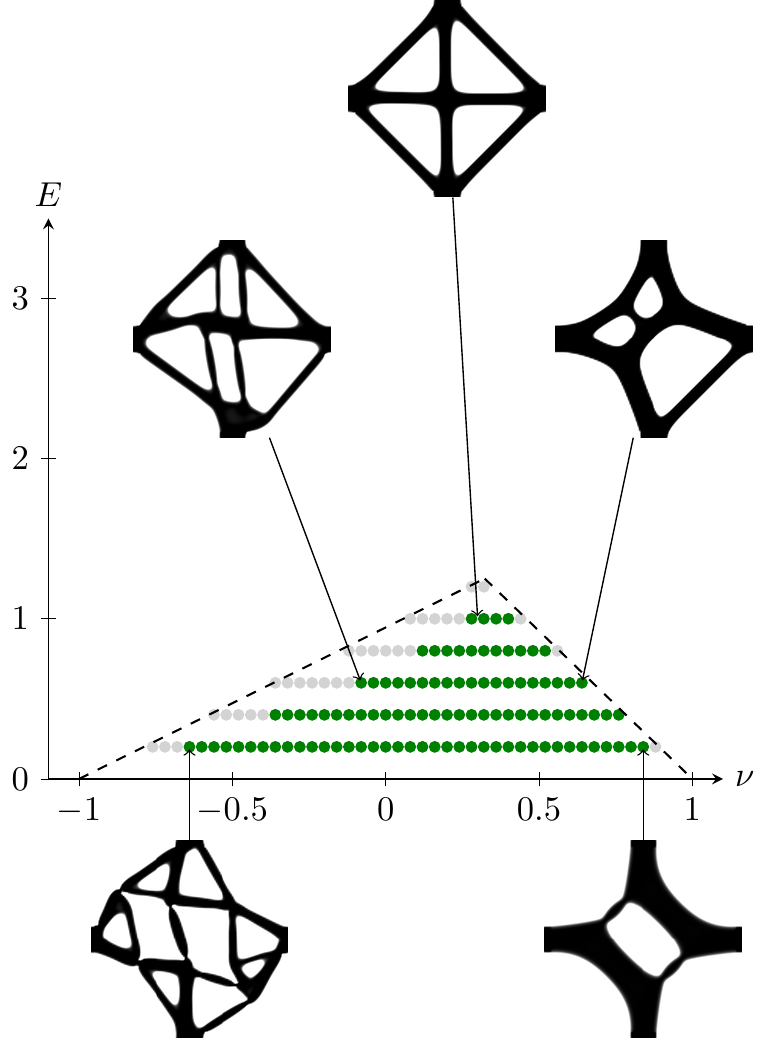}
 \includegraphics[width=0.245\textwidth]{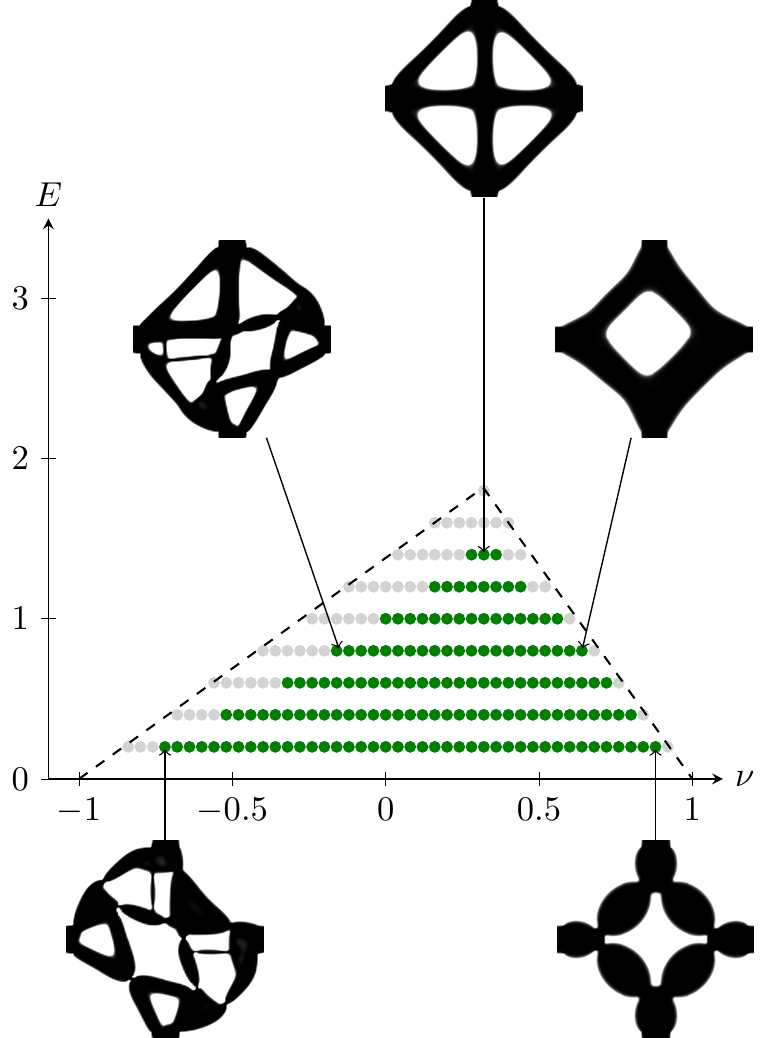}
 \includegraphics[width=0.245\textwidth]{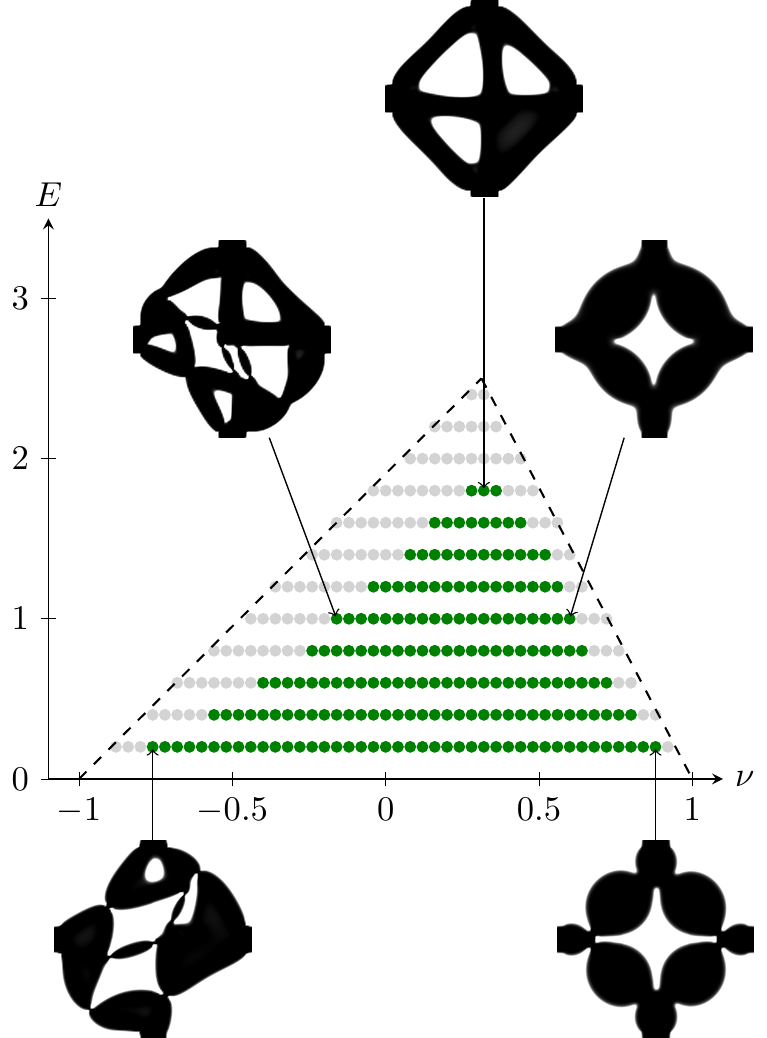}
 \includegraphics[width=0.245\textwidth]{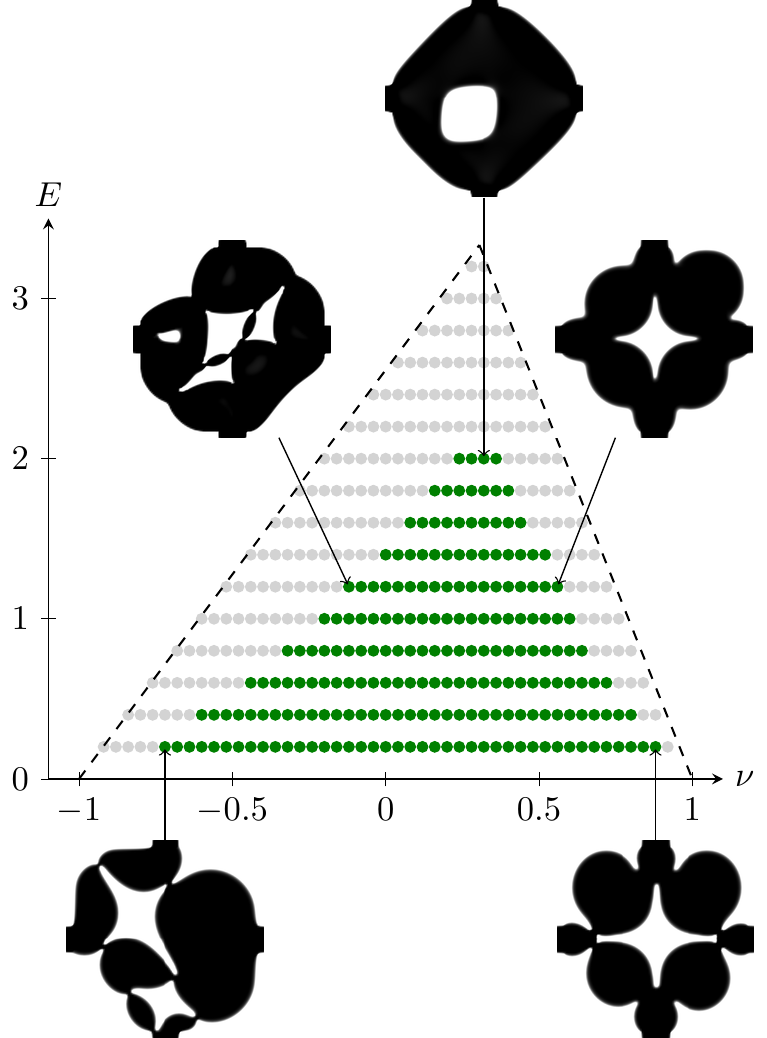}
}
\caption{
Set of realizable effective elasticity tensors (green point clouds) 
for different fixed volume fractions 
$\theta_s = 0.3$, $0.4$, $0.5$ and $0.6$ (from left to right) of the stiff material.
Dashed lines correspond to the upper Hashin-Shtrikman bounds for the specific $\theta_s$.
For points colored in grey inside these bounds, we computationally did not observe an admissible fine structure with our optimization procedure.
Selected fine-scale structures are also shown for some extremal points. 
}
 \label{fig:CellFixedVol}
\end{figure}

\paragraph*{Optimization for fixed volume fraction in 2d.}
First, we consider the optimization problem \eqref{eq:shapeOptProblemCell}
for a parameter $c_V=0$, but with an additional volume constraint on the stiff phase. 
Note that we later do not use this specific type of problem.
We rather depict these solutions for comparison, 
since the classical G-closure set in the literature~\cite{Ch00a} is usually considered for a fixed volume fraction.
In Figure~\ref{fig:CellFixedVol}, we show the resulting point clouds in a $(\nu,E)$-coordinate system for specific volume fractions $\theta_s$ of the stiff phase.

\begin{figure}[htbp]
\centering{
 \includegraphics[width=0.29\textwidth]{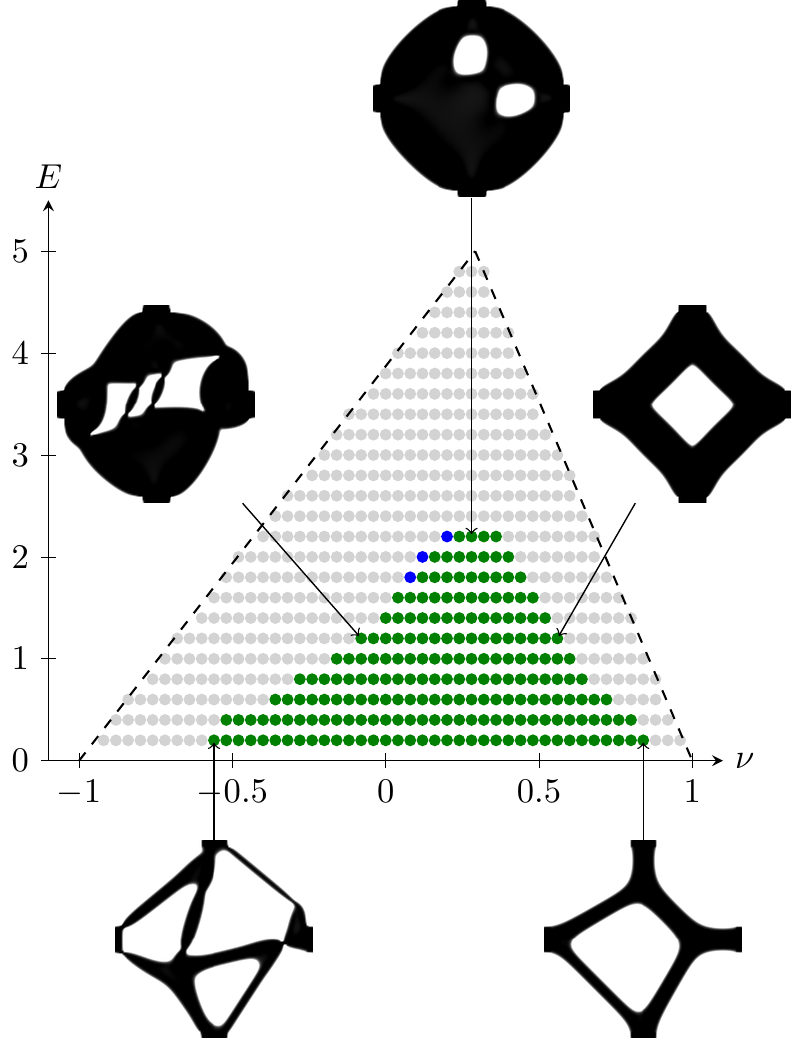}
 \hspace*{0.05\textwidth}
 \includegraphics[width=0.29\textwidth]{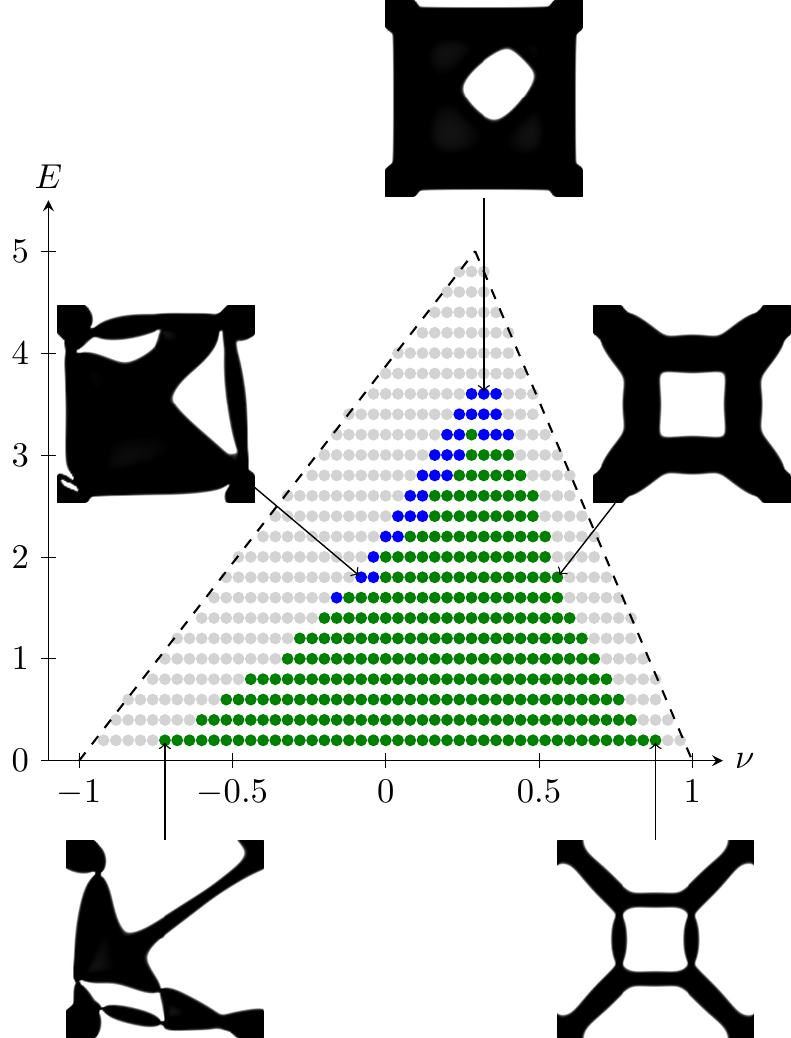}
 \hspace*{0.05\textwidth}
 \includegraphics[width=0.29\textwidth]{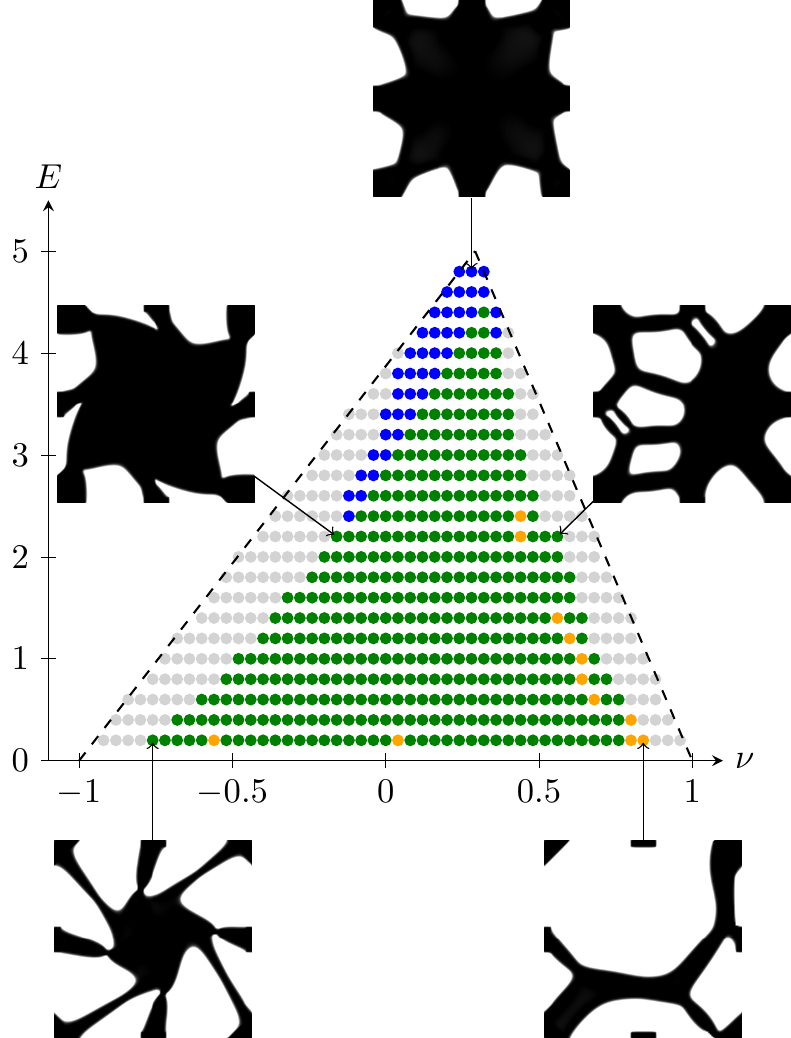}
 }
 \caption{
   We compare point clouds approximating the set of realizable manufacturable materials (green) for several manufacturability constraints.
   Certain homogenized material tensors are manufacturable, but require a higher volume fraction than the initially chosen $\theta_s = 0.75$ (blue).
   In the case of bridges both in the corners and the center of the edges (right), we obtain certain disconnected structures for the hard phase (yellow) which do not use all bridges.
   Trial points inside the Hashin--Shtrikman bounds with no admissible solution are colored in gray.
 }
 \label{fig:CellOptVol}
\end{figure}

\paragraph*{Optimization of the volume fraction in 2d.}
Now, we consider the optimization problem \eqref{eq:shapeOptProblemCell} for a parameter $c_V=1$, 
\ie we optimize the volume fraction of the stiff phase and add the perimeter functional to regularize the interface
with a penalty factor $c_P = 0.05$.
In Figure~\ref{fig:CellOptVol} we compare three different manufacturability constraints given by different bridge sets $B^1$ and $B^2$.
For the upper Hashin--Shtrikman bounds, which we again show as dashed lines, we choose a volume fraction of $\theta_s = 0.75$.

\begin{figure}[htbp]
\centering{
 \includegraphics[width=0.5\textwidth]{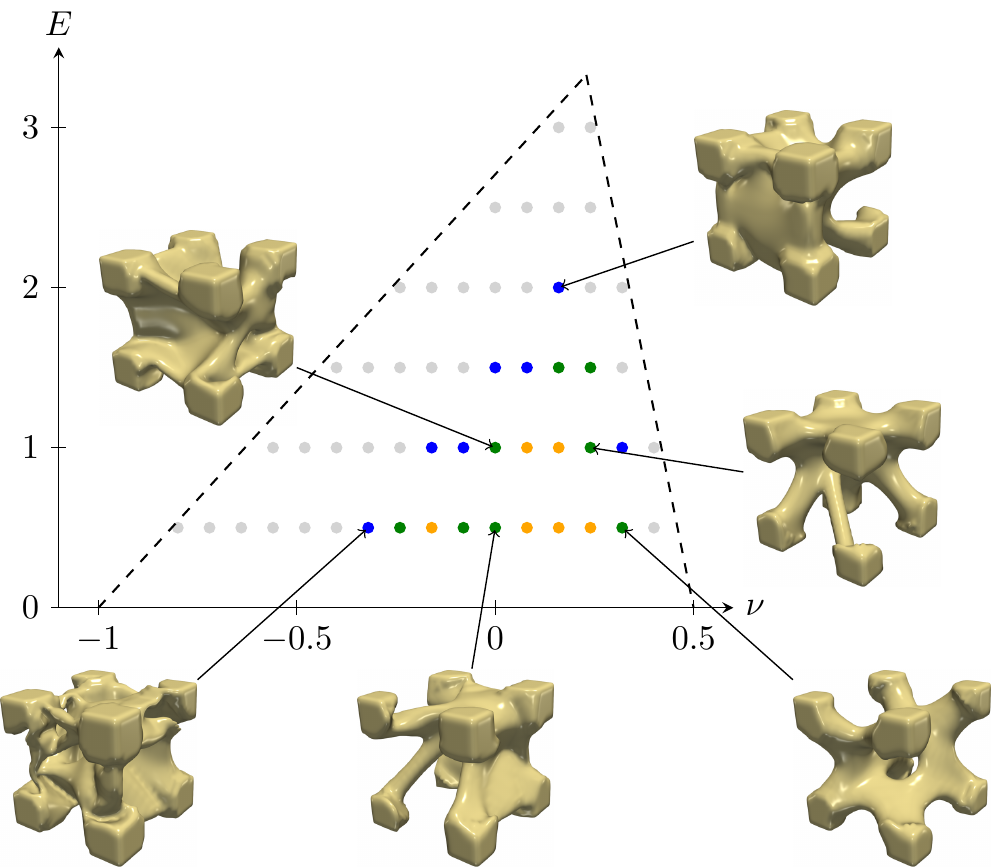}
 }
 \caption{
   We depict a point cloud in the $(\nu,E)$-coordinate system
   approximating the set of realizable manufacturable materials in 3d together with selected images of the associated hard phase.
 }
 \label{fig:CellOptVol3d}
\end{figure}

\paragraph*{Optimization of the volume fraction in 3d.}
For the discretization of the unit cube $[0,1]^3 \subset \R^3$, we use a mesh size of $h=2^{-5}$.
As for the other computations in 2d, 
we choose a Poisson's ratio $\nu = 0.25$, a Young's modulus $E = 10$,
and $\delta = 10^{-4}$ as scaling factor to approximate the void material.
For the material bridges, we choose the set $B^1$ at the eight corners of the unit cube.
In Figure~\ref{fig:CellOptVol3d}, we show the resulting point cloud in the $(\nu,E)$-coordinate system.
As in Figure~\ref{fig:CellOptVol}, certain homogenized material tensors are manufacturable, but require a higher volume fraction than the initially chosen $\theta_s = 0.5$ (blue).
Furthermore, we obtain certain disconnected hard phase structures (yellow).
Compared to the two-dimensional case, we remark that the corresponding point cloud for a bridge set $B^1$ in the six midfaces is essentially smaller, 
whereas using bridges both around the corners and around the midpoints of the faces leads to a point cloud with essentially disconnected structures.
By a suitable reinitialization of initially disconnected, optimal microstructures (yellow dots) we will be able to enlarge the range of manufacturable material parameters (\cf~Figure~\ref{fig:samplingexample} in the next section).

\section{Restricted free material optimization on the macro-scale} \label{sec:macroopt}
In this section, we consider an optimization problem on the macro-scale, where we optimize the 
macroscopic, linearized elasticity model in a 
restricted
free material optimization sense over elasticity tensor 
distributions $\Etensor \colon \domain \to \R^{d^4};\ x\mapsto \Etensor[x]$ with $\Etensor[x]$ 
in a subclass of effective elasticity tensors arising from microscopic material patterns 
as described in Section~\ref{sec:twoscale} above. 
Based on computational results of the offline microscopic optimization this free material optimization is performed very efficiently online for given load scenarios.
To this end, we consider a macroscopic displacement $\displaceMacro \in \displaceMacro^\partial + H^{1,2}_{\Gamma_D} (\domain)^{d}$ 
and a macroscopic total free energy 
$\Etot[\Etensor,\displaceMacro] \coloneqq \Eelast[\Etensor,\displaceMacro] - \Epot[\displaceMacro]$  
involving the stored elastic energy and the potential energy 
\begin{align}
 \Eelast[\Etensor,\displaceMacro] \coloneqq  \tfrac12 \int_\domain \Etensor \varepsilon[\displaceMacro] : \varepsilon[\displaceMacro] \d x  \,, \quad
 \Epot[\displaceMacro] \coloneqq \int_\domain f \cdot \displaceMacro \d x + \int_{\Gamma_N}  g \cdot \displaceMacro \d a
\end{align}
as in Section~\ref{sec:twoscale} above where now
\begin{align}
 \Etensor(x) \in \mathcal{C}_{\text{real}} \coloneqq \{\Etensor \in \mathcal{C}_{\text{adm}} \ \vert \ \Etensor = \Etensor^\ast[\chi]  \,, \ \chi \in \mathcal{B} \}
\end{align}
and $\mathcal{C}_{\text{adm}}$ is the set of admissible elasticity tensors. 
As we have mentioned in Section~\ref{sec:micro-scale}, 
in our application, we will consider for simplicity $\mathcal{C}_{\text{adm}}$ as the (two-dimensional) space of isotropic elasticity tensors, 
since the full space of all possible elasticity tensors is too high dimensional. 
For this given $\mathcal{C}_{\text{adm}}$,
the set $\mathcal{C}_{\text{real}}$ has been identified in the offline phase of our optimization approach as explained above.
So far, note that we have only computed a discrete point cloud of the set $\mathcal{C}_{\text{real}}$,
\st a suitable interpolation is required, which we will detail later in this section.
We denote by $\displaceMacro[\Etensor]$ the minimizer of the total elastic energy over all 
$\displaceMacro \in \displaceMacro^\partial + H^{1,2}_{\Gamma_D} (\domain)^d$.
This couples the macroscopic free material optimization to the microscopic scale,
since $\displaceMacro[\Etensor]$ is the weak solution of the linearized elasticity problem with the 
inhomogeneous elasticity tensor field $\Etensor$ point-wise a.e. in $\mathcal{C}_{\text{real}}$, 
which means that the elasticity tensor can be realized as an effective elasticity tensor on the micro-scale.

Taking into account a cost functional $\mathcal{J}[\displaceMacro]$,
we ask for a minimizer of $\mathcal{J}$ subject to the constraint that $\displaceMacro = \displaceMacro[\Etensor]$.
Possible macroscopic cost functionals are of compliance or tracking type, choosing 
\begin{align}
 \mathcal{J}_{\text{compl}}[\displaceMacro] \coloneqq \Epot[\displaceMacro] \quad \text{or } \quad
 \mathcal{J}_{\text{track}}[\displaceMacro] \coloneqq \int_{\domain^{\text{track}}} |\displaceMacro - \displaceMacro^{0}|^2 \d x
\end{align}
for a given displacement $\displaceMacro^{0}$ on a tracking domain $\domain^{\text{track}}\subset \domain$.
In addition, we might want to minimize the weight of the material phase with characteristic function $\chi$
or the density of the microscopic perimeter of the interface between material phases, or a linear combination of the two.
This is done via the function $\mathbf{j}[\Etensor]$,
which was defined 
in \eqref{eq:localCost} 
as the sum (with suitable coefficients $c_V$ and $\widehat{c_P}$) 
of the microscopic material density and the microscopic perimeter density corresponding to the considered
microscopically optimal material distribution $\chi \in \mathcal{B}$ 
for which $\Etensor$ is the effective elasticity tensor.
For details we refer to Section~\ref{sec:micro-scale} on the microscopic optimization problem.
The resulting total cost functional is then given by
\begin{align}
 \mathbf{J}_{\text{macro}}[\Etensor] 
 = \int_\domain \mathbf{j}[\Etensor(x)] \d x + \mathcal{J}_{\text{mech}}[\displaceMacro[\Etensor]]
\end{align}
with the mechanical cost functional $\mathcal{J}_{\text{mech}}$ being either $\mathcal{J}_{\text{compl}}$ or $\mathcal{J}_{\text{track}}$
and $\mathbf{j}: \R^{d^4} \to \R^+_0\,,\ \Etensor \mapsto \mathbf{j}[\Etensor]$ 
being the cost density which reflects the cost functional evaluated on the micro-scale and has already been computed in the offline phase (\cf~Section~\ref{sec:micro-scale}).
We recall that $\mathbf{j}[\Etensor] = \infty$
for $\Etensor \notin \mathcal{C}_{\text{real}}$.

Alternatively, for the macroscopic optimization problem,
we might take into account the equality constraint on the integrated cost density
$\int_\domain \mathbf{j}[\Etensor(x)] \d x = \mathbf{j}_{\text{constr}}$
and consider for the total cost functional just 
the compliance functional $\mathcal{J}_{\text{compl}}$ 
or the tracking functional $\mathcal{J}_{\text{track}}$.

\paragraph*{Approximate fine-scale realization.}
In what follows we show how to design a manufacturable workpiece with a fine-scale geometry which reflects 
the free material optimization on the macro-scale with a piecewise constant microscopically realizable
effective elasticity tensor on a selected fine-scale. 
To this end, we proceed as follows.
We assume that the domain $\domain$ can be decomposed as 
\begin{align}\label{eq:decompositionMacroDomain}
 \domain =  \bigcup_{\alpha\in I_\epsilon} \epsilon (\alpha+Q) \, ,
\end{align}
where $I_\epsilon$ is a set of multi-indices in $\mathbb{Z}^d$ corresponding to a fixed fine-scale parameter $\epsilon$.
We consider piecewise constant elasticity tensors 
$\Etensor:I_\epsilon\to\R^{d^4}$, which can be identified with elements of $L^\infty(\domain)^{d^4}$ which are  constant on each cell $\epsilon (\alpha+Q)$.
Then, given the field of optimal, effective elasticity tensors $\Etensor^\ast$,
we construct a material pattern on the cell $\epsilon (\alpha+Q)$ with associated characteristic function 
\begin{align}
 \chi^{\text{macro}}(x) = \chi^\ast[\Etensor^\ast(\alpha)]\left(\tfrac{x}{\epsilon} - \alpha\right)
\end{align}
for all $x\in \epsilon (\alpha+Q)$.
Here, $\chi^\ast[\Etensor^\ast(\alpha)]$ is the characteristic function solving the 
microscopic shape optimization problem with effective elasticity tensor $\Etensor^\ast(\alpha)$.
The printability constraint introduced above ensures that each component of the bridge set close to the boundary of a cell $\epsilon (\alpha+Q)$
has a counterpart in the neighbouring cells $\epsilon (\alpha \pm e_k +Q)$ 
for $k\in\{1,\ldots, d\}$ (as long as $\epsilon (\alpha \pm e_k +Q) \subset D$)
independently of the selected elasticity tensor on the cell.
Let us remark that the selection of a fixed scale $\epsilon$ to define a manufacturable geometry on that scale 
conceptually differs from the periodicity requirement in the scaling limit of classical homogenization theory in the second, fast varying variable of $(x,y) \mapsto \chi(x,y)$. 

Taking into account the volume $\epsilon^d$ of the cells $\epsilon (\alpha+Q)$
and the perimeter scale factor $\epsilon^{d-1}$ of the characteristic function $\chi^{\text{macro}}$
on $\tfrac{\alpha+Q}{\epsilon}$ 
compared to the perimeter of $\chi^\ast[\Etensor^\ast(\alpha)]$ on 
$Q$ we observe that 
\begin{align}
\sum_{\alpha\in I_\epsilon} \epsilon^d\mathbf{j}[\Etensor^\ast(\alpha)]=
 \int_\domain \mathbf{j}[\Etensor^\ast(x)] \d x = c_V \int_D \chi^{\text{macro}}(x) \d x + \widehat{c_P} \, \epsilon\, |D \chi^{\text{macro}}|(\domain)  \,.
\end{align}
Hence, for fine-scale cells $\epsilon (\alpha+Q)$ of fixed scale $\epsilon$ the cost component
$\int_\domain \mathbf{j}[\Etensor(x)] \d x$ is indeed the sum of the total volume of the hard material phase in $D$ 
with weight $c_V$ and the total perimeter of this phase with weight $\widehat{c_P} \epsilon$.

\paragraph*{B-Spline interpolation of the set of realizable elasticity parameters.}
In the optimization on the macro-scale, we currently select elasticity tensors $\Etensor$ from the set of 
realizable elasticity tensors $\mathcal{C}_{\text{real}}$.
The associated cost is given by the cost density function $\mathbf{j}[\cdot]$. 
Neither the set $\mathcal{C}_{\text{real}}$ nor the microscopic functional $\mathbf{j}[\cdot]$ can be given explicitly. 
In fact, on a larger set of elasticity tensors $\Etensor$, we optimize over phase-field functions $\pf$
subject to the constraint that $\Etensor$ is the associated effective elasticity tensor 
for the material distribution on the fundamental cell $Q$ on the micro-scale described by $\pf$, 
\ie $\Etensor=\Etensor^\ast[\pf]$. 
If there is no such phase-field, then $\Etensor$ belongs to the complement of $\mathcal{C}_{\text{real}}$.
Having identified a phase-field $\pf$, the cost density $\mathbf{j}[\Etensor]$ can be computed. 
Thus, we are only able to generate samples $(\Etensor_i, \mathbf{j})$ for $i$ in some sample index set and 
with $\mathbf{j}[\Etensor_i] = \infty$ for $\Etensor_i \notin \mathcal{C}_{\text{real}}$.
Since so far we only have a discrete set of samples, 
we require a suitable smooth interpolation strategy
to appropriately perform the above optimization. 

\begin{figure}[htbp]
\includegraphics[width=1.\textwidth]{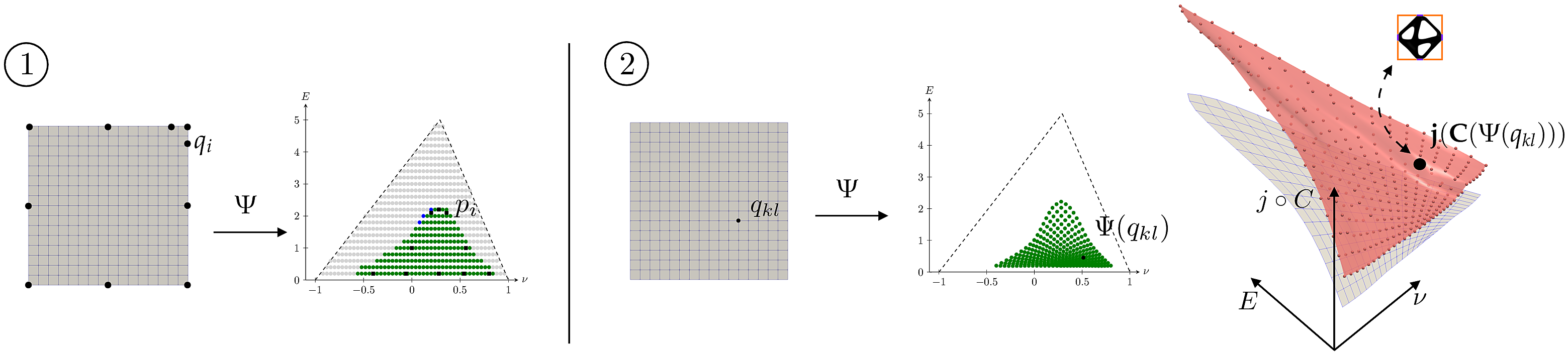}
\caption{
A sketch of resampling and B-spline interpolation approach:
on the left step $1$ is depicted with the selection of sample points $p_i$ 
close to the boundary $\partial \mathcal{P}_{\text{real}}$ and their preimages $q_i$ on $[0,1]^2$ which give rise to constraints in the
variational definition of the B-spline parametrization $\Psi$ of an interior approximation of $\mathcal{P}_{\text{real}}$;
on the right step $2$ is illustrated with lattice points $q_{kl}$ on $[0,1]^2$ and their image points $\Psi(q_{kl})$ as well as the associated 
cost values  $\mathbf{j}(\Etensor(\Psi(q_{kl})))$ above the points $\Psi(q_{kl})$ in the $(\nu,E)$ plane. These cost values are associated with an optimal microscopic material pattern.  
} 
\label{fig:overview}
\end{figure}

We assume that the set of admissible elasticity tensors
$\mathcal{C}_{\text{adm}}$ can be parametrized over a small set of parameters $p$ in a subset $\mathcal{P}_{\text{adm}}$
of admissible parameters in a parameter space $\mathcal{P}$. 
In our application, we consider isotropic elasticity tensors as admissible and parametrize them as $p \mapsto \Etensor(p)$ according to  
\eqref{eq:isotropC} with $p=(\nu, E)$. 
In turn, the set of realizable elasticity tensors $\mathcal{C}_{\text{real}}$,
which arise from a material
pattern represented by some $\chi \in \mathcal{B}$,
corresponds to a specific set of 
realizable effective 
parameters, 
$\mathcal{P}_{\text{real}}\subseteq
\mathcal{P}_{\text{adm}}$, with
$\mathcal{C}_{\text{real}} = \Etensor(\mathcal{P}_{\text{real}})\subseteq \mathcal{C}_{\text{adm}}$
(\cf the numerical results in Section~\ref{sec:fineScaleResults} or  Figure~\ref{fig:overview}).
Due to the Hashin--Shtrikman bounds~\cite{HaSh63}, the search for $\mathcal{P}_{\text{real}}$ can be 
restricted to the values of  $p=(\nu, E)$ contained in a certain triangle. 
The various constraints we introduced, including for example the bridges on the boundary, will ensure that $\mathcal{P}_{\text{real}}$ is a strict subset of this triangle.

\medskip

\begin{figure}[htbp]
 \centering{
   \includegraphics[width=0.2\textwidth]{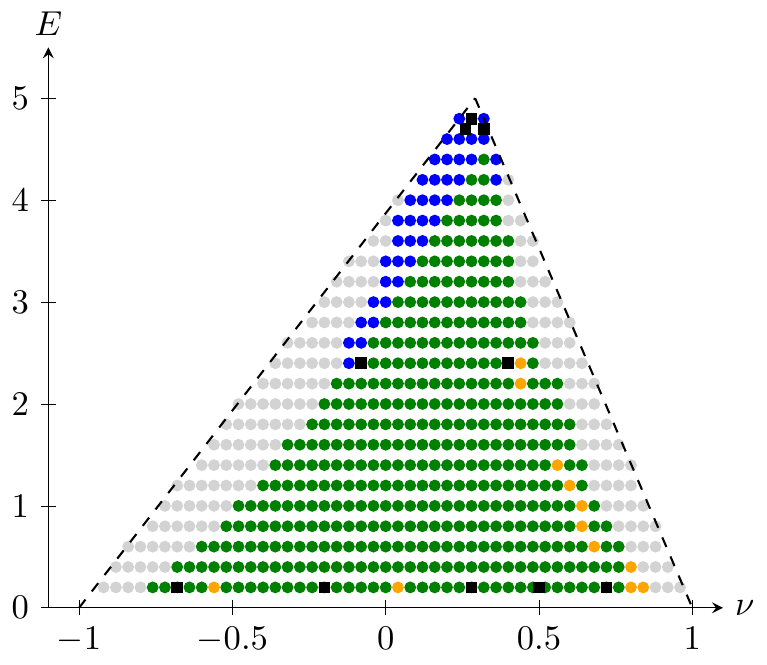}
   \includegraphics[width=0.2\textwidth]{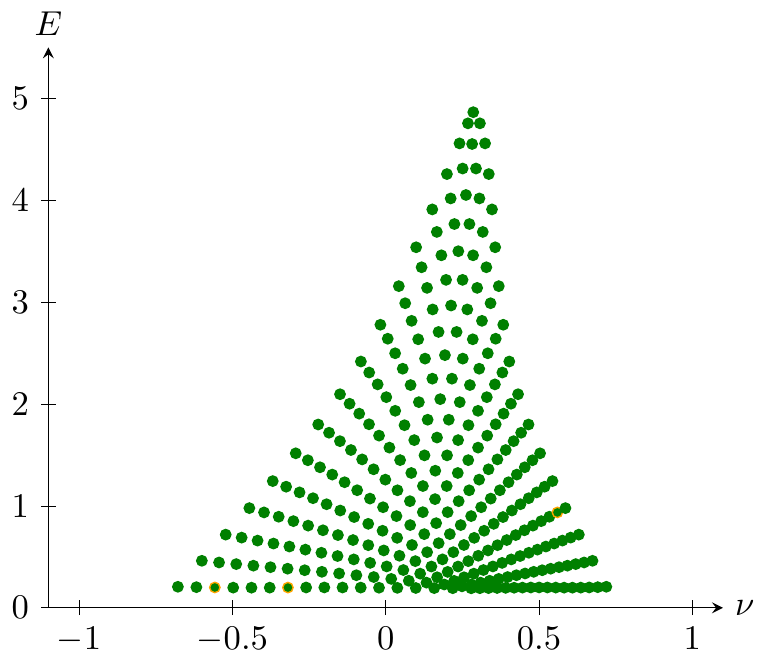}
   \hspace*{0.1\textwidth}
   \includegraphics[width=0.2\textwidth]{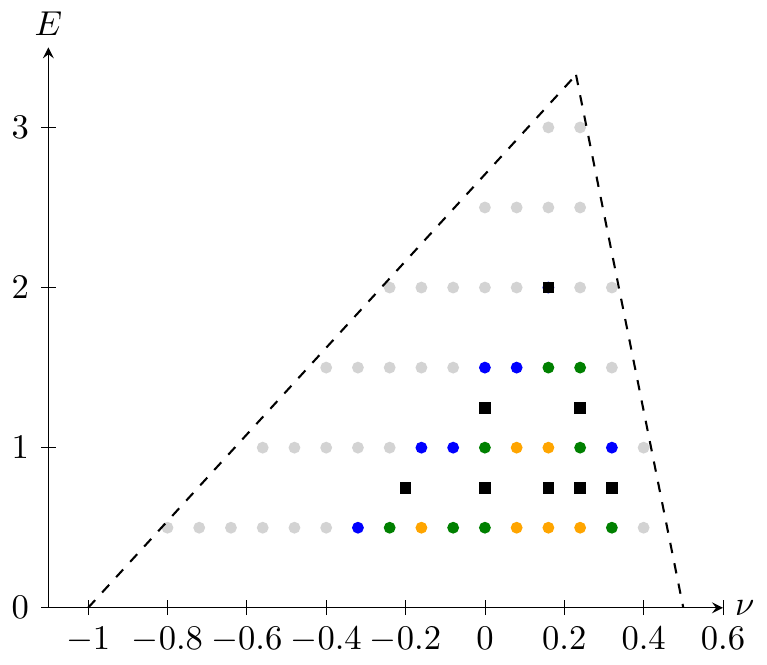}
   \includegraphics[width=0.2\textwidth]{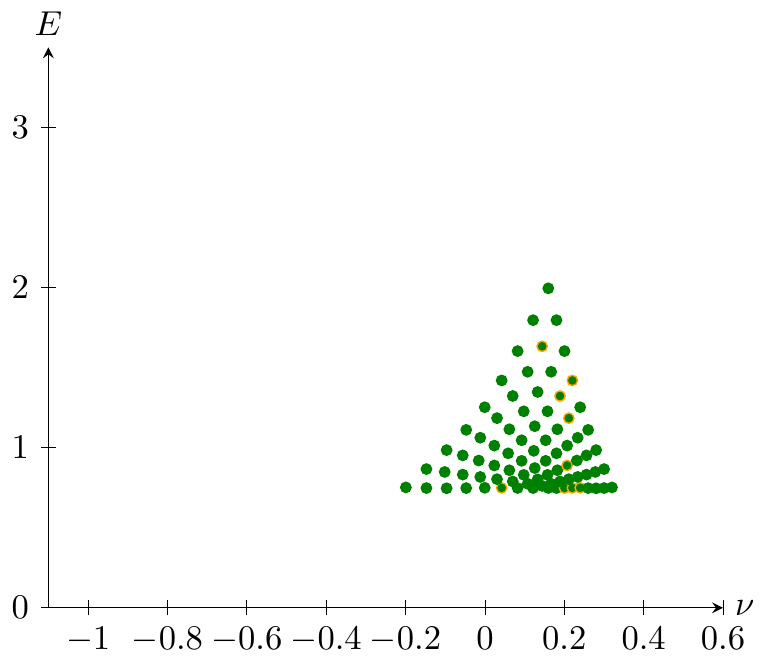}
}
\caption{
We show the resampling of point clouds for 2D bridge sets at vertices and midfaces (left) and for 3D bridge sets at vertices (right).
Here, we indicate sampling points by black boxes, disconnected optimized hard phase microstructures by yellow dots, 
The material parameters of previously disconnected structures could be recovered by connected structures (green dots
on top of yellow dots) via initialization of the phase-fields in the optimization approach with convex combinations of already computed microstructures which are connected. 
}
\label{fig:samplingexample}
\end{figure}

\noindent Computationally, we now proceed in two steps (\cf Figures~\ref{fig:overview} and \ref{fig:samplingexample}):
\medskip

\begin{enumerate}[label=B\arabic*,ref=B\arabic*]
 \item \label{computationBSplineStepOne}
  In a first step, we compute a smooth approximation of $\mathcal{P}_{\text{real}}$ via 
  a cubic B-spline parametrization over the unit square. 
  More precisely, we discretize the unit interval $[0,1]$ uniformly with mesh size $\gridSizeSplineChart$,
  consider associated cubic B-splines and denote by 
  $\splineFESpace$ the space of tensor products of these cubic B-splines on the unit interval.
  Note that a function in $\splineFESpace$ has global $C^1$-regularity on the unit square $[0,1]^2$.
  Now, we select by hand a set of sample points $p_i$  
  close to the boundary $\partial \mathcal{P}_{\text{real}}$ and preimages $q_i$ on $[0,1]^2$ 
  for $i$ in some index set $I^\text{sample}$ (see Figure~\ref{fig:samplingexample}).
  Here, we typically chose 
  $p(1,0)$ with large Poisson's ratio,
  $p(0,1)$ with small Poisson's ratio,
  $p(1,1)$ with large Young's modulus,
  and $p(0,0)$ with intermediate Poisson's ratio and small Young's modulus.
  Even if we choose further points $p_i$ (\cf Figures~\ref{fig:overview}), we did not intend to explore the geometry of the boundary $\partial \mathcal{P}_{\text{real}}$.
  Next, we define an interpolating deformation $\Psi$ on $[0,1]^2$ as
  \begin{align}
  \Psi = \argmin_{\widetilde \Psi \in \splineFESpace^2 \; : \; \widetilde\Psi(q_i) = p_i \, \forall i\in I^\text{sample}} \int_{[0,1]^2} |D^2 \widetilde \Psi|^2 \d p \, .
  \end{align}
  In the implementation, the integral is approximated by a tensor product Gauss quadrature with $15^2$ quadrature points per element.
  With the choice of the samples $p_i$ and the corresponding $q_i$, 
  we ensure that $\Psi([0,1]^2)$ is a suitable subset of the set $\mathcal{P}_{\text{real}}$
  obtained numerically via the finite element discretization of 
  our phase-field approach.

 \item \label{computationBSplineStepTwo}
  In the next step, we resample the mapping from elasticity tensors to cost densities.
  To this end, we consider the regular lattice of B-spline control points $q_{kl}$ on $[0,1]^2$ with $k,\ l \in \{1,\ldots, N\}$
  corresponding to the mesh size $\gridSizeSplineChart$.
  For each point $q_{kl}$, we compute the cost density $\mathbf{j}(\Etensor(\Psi(q_{kl})))$ numerically.
  Finally, we define the B-spline interpolation 
  $\mathbf{j}_{\text{ref}} \coloneqq \mathbf{j}\circ \Etensor\circ \Psi$ 
  of these densities interpolating the values $\mathbf{j}(\Etensor(\Psi(q_{kl})))$ at the lattice points $q_{kl}$ on the reference domain $[0,1]^2$.
\end{enumerate}

\paragraph{Reformulation of the macroscopic free material optimization problem.}
For the free material optimization problem on the macro-scale, the cost functional can be rewritten as 
\begin{align}
 \widehat{\mathbf{J}}_{\text{macro}}[\mathbf{q}] 
 \coloneqq \int_\domain \mathbf{j}_{\text{ref}}(\mathbf{q}(x)) \d x 
         + \mathcal{J}_\text{mech}\left[\displaceMacro[\Etensor(\Psi \circ \mathbf{q})]\right] \, ,
\end{align}
which we aim at minimizing over the map $\mathbf{q} \colon \domain \to [0,1]^2;\ x \mapsto \mathbf{q}(x)$.
Here, remember that $\mathbf{j}_{\text{ref}}$ is the B-spline interpolation as defined in \ref{computationBSplineStepTwo}.
The highly nonlinear constraint that the elasticity tensor $\Etensor$ has to be realizable with finite microscopic cost density 
is now replaced by the simple box constraints $\mathbf{q}(x) \in [0,1]^2$ for all $x \in \domain$, 
which is much easier to handle numerically in an interior point approach for the free material optimization on the macro-scale. 
Finally, we remark that starting from scratch with random initial phase-fields we sometimes observe optimized hard phase microstructures which are disconnected.
These usually leave one of the bridge sets in $B^1$ (\cf Figure~\ref{fig:CellOptVol} in 2D and \cf Figure~\ref{fig:CellOptVol3d} in 3D) disconnected from the remaining structure.
Frequently, in the B-spline resampling step, connected optimal microstructures can be computed for these parameter sets by using suitable convex combinations of neighbouring sample points as initialization.
More precisely, for the phase-field, we take into account an affine interpolation of the phase-field on neighbouring sample points for which optimal connected microstructure were already obtained (\cf~Figure~\ref{fig:samplingexample}).
 
\paragraph*{Numerical implementation of the macroscopic free material optimization problem.}
Finally, we have to discretize the macroscopic free material optimization problem.
Since we made in \eqref{eq:decompositionMacroDomain} the assumption
that our macroscopic domain $\domain$ can be decomposed into scaled unit cubes,
we use, as for the fine-scale optimization problem, 
a quadrilateral mesh in the case $d=2$
and a hexahedral mesh in the case $d=3$.
Note that the macroscopic mesh size is given by $\macrogridsize=\epsilon$, s.t. the corresponding elements are given by the fine-scale cells $\macrogridsize(\alpha + Q)$.
Then we define the space $\mathcal{V}_\macrogridsize$ of piecewise $d$-affine, continuous functions 
and consider this space as the ansatz space for the discrete macroscopic displacements.
For the assembly of the stiffness matrix and the right hand side we apply a tensor product Simpson quadrature rule.
Finally, the actual shape optimization problem in the unknown given by a piecewise constant $\mathbf{q}$ is solved using the IPOPT package \cite{WaBi06}.

\section{Numerical results for the two-scale shape optimization problem} \label{sec:twoScaleResults}
In this section, we present numerical results for the macroscopic shape optimization problem taking into account the admissible set of manufacturable microstructures.
Remember that this admissible set was computed in Section~\ref{sec:fineScaleResults} on the unit cube $[0,1]^d$ with a certain mesh size $h$ ($h=2^{-7}$ in 2d and $h=2^{-5}$ in 3d).
For the material to be manufacturable, we have chosen $E = 10$ for Young's modulus, $\nu = 0.25$ for the Poisson's ratio, and $\delta = 10^{-4}$ as factor to approximate the void material.
Moreover, for the spline finite element functions $\Psi$
as described in step \ref{computationBSplineStepOne},
we use a mesh size $\gridSizeSplineChart = 2^4$ for $d=2$ and $\gridSizeSplineChart = 2^3$ for $d=3$.
In the visualization, for an interpolated effective property, we have just depicted the microstructure of the closest point in the sampling.
  
\begin{figure}[htbp]
 {\includegraphics[width=\textwidth]{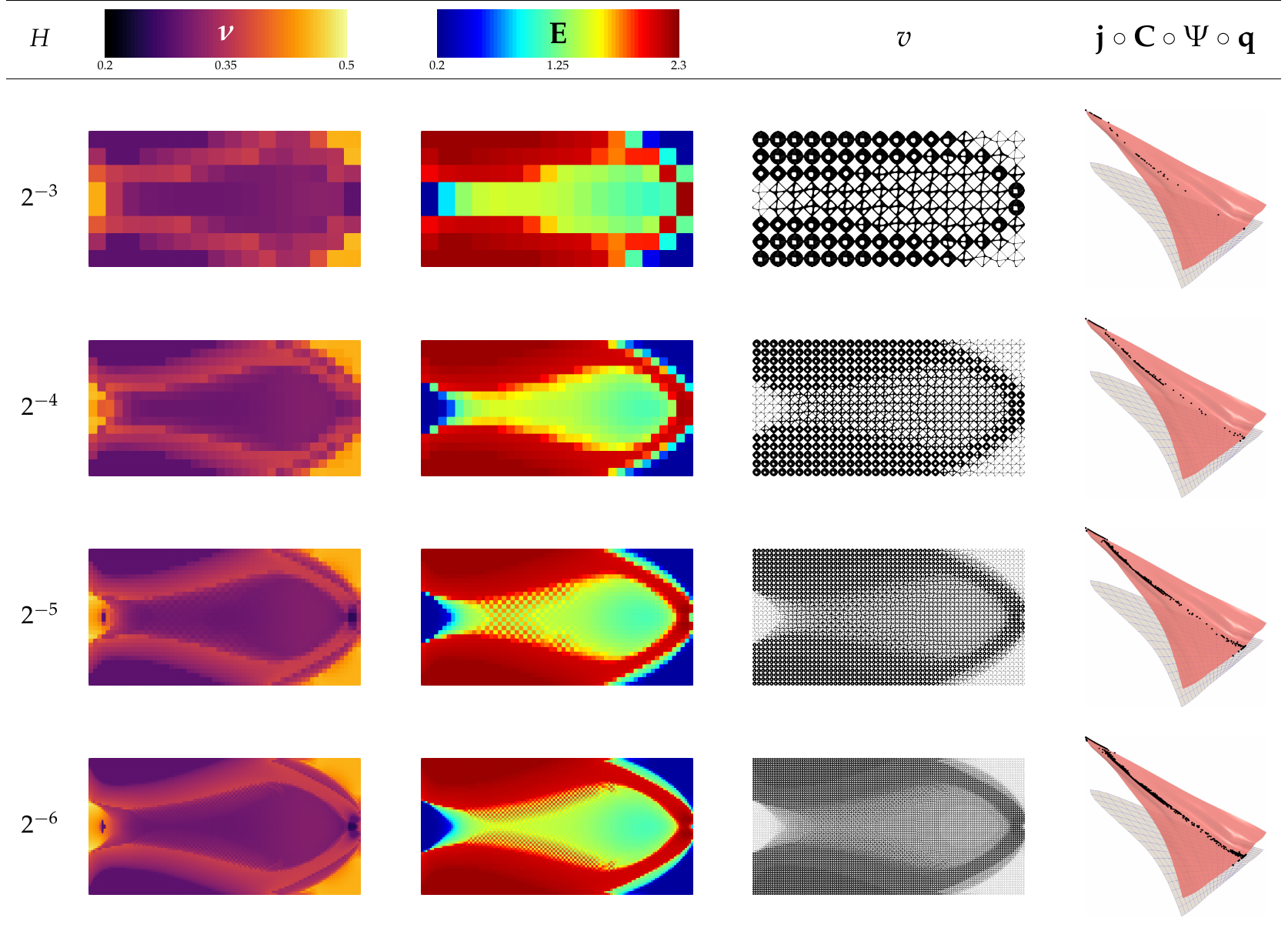}} 
 
 \vspace*{2ex}
 
 \centering{
    \begin{minipage}{0.45\textwidth}
        {\includegraphics[width=0.95\textwidth]{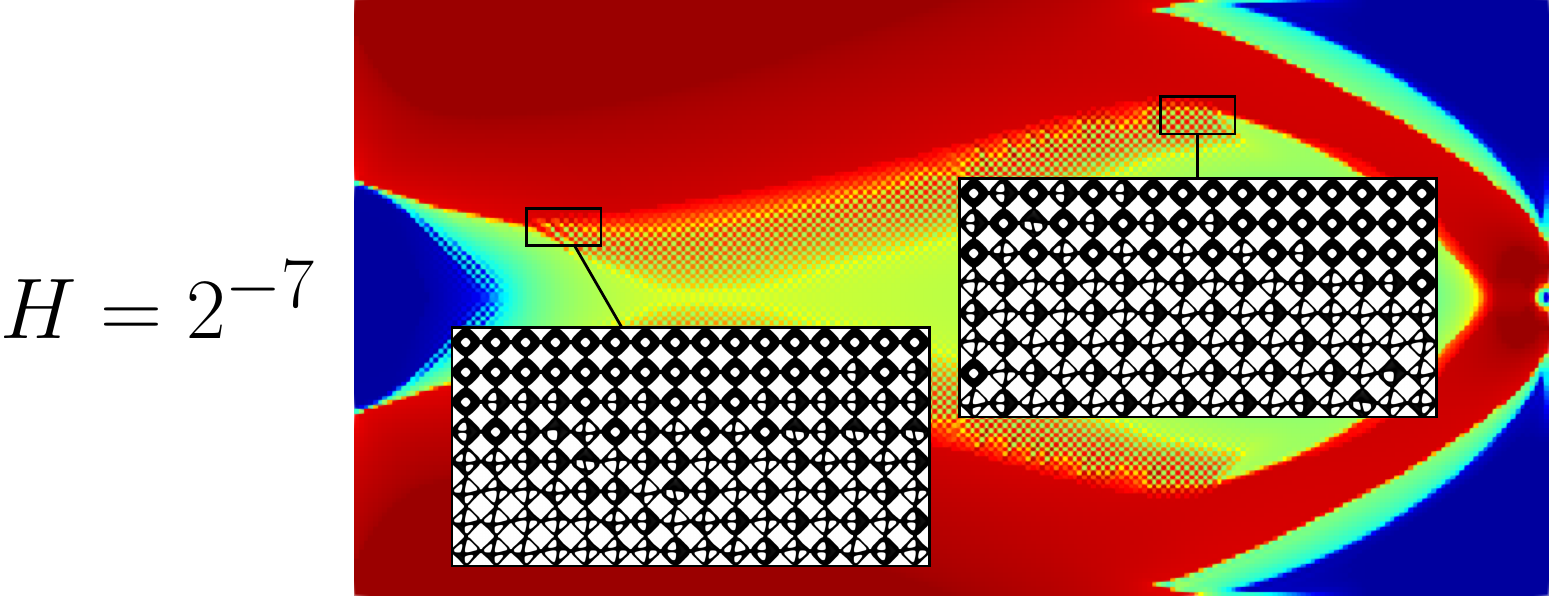}} 
    \end{minipage}
    \begin{minipage}{0.45\textwidth}
        {\includegraphics[width=0.95\textwidth]{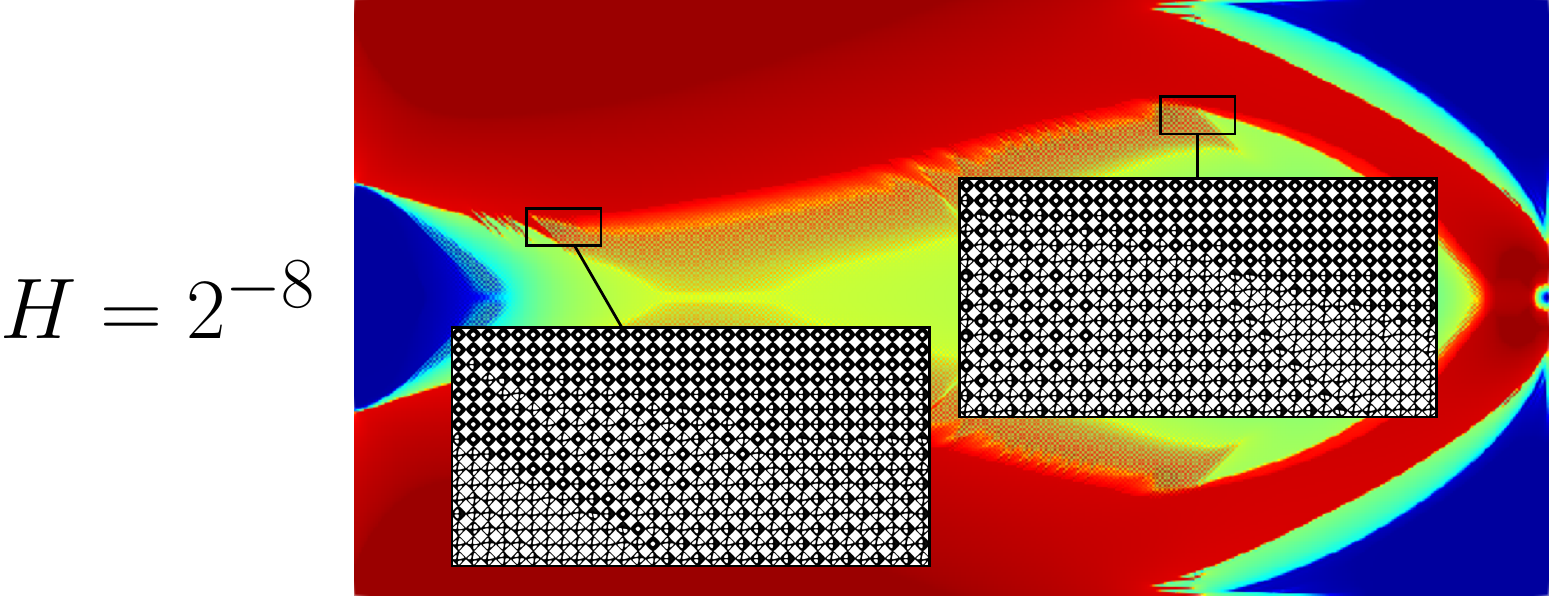}} 
    \end{minipage}
 }
\caption[2D Cantilever, bridges mid, vol 1]{
Computational results are shown for the cantilever problem for a volume constraint $\Volume_H = 1$.
For different mesh sizes, 
we depict Poisson's ratio and Young's modulus on the macro-scale, 
the full two-scale structure,
and the cost values of the used microcells.
Here, the admissible fine-scale structures have material bridges at the midfaces. 
In the bottom row, we explore local microstructures for macroscopic mesh sizes $H=2^{-7}$ (left) and $H=2^{-8}$ (right) 
at selected regions.
}
\label{fig:TwoScaleCantilverBridgesMiddleVol1}
\end{figure}

\begin{figure}[htbp]
 {\includegraphics[width=\textwidth]{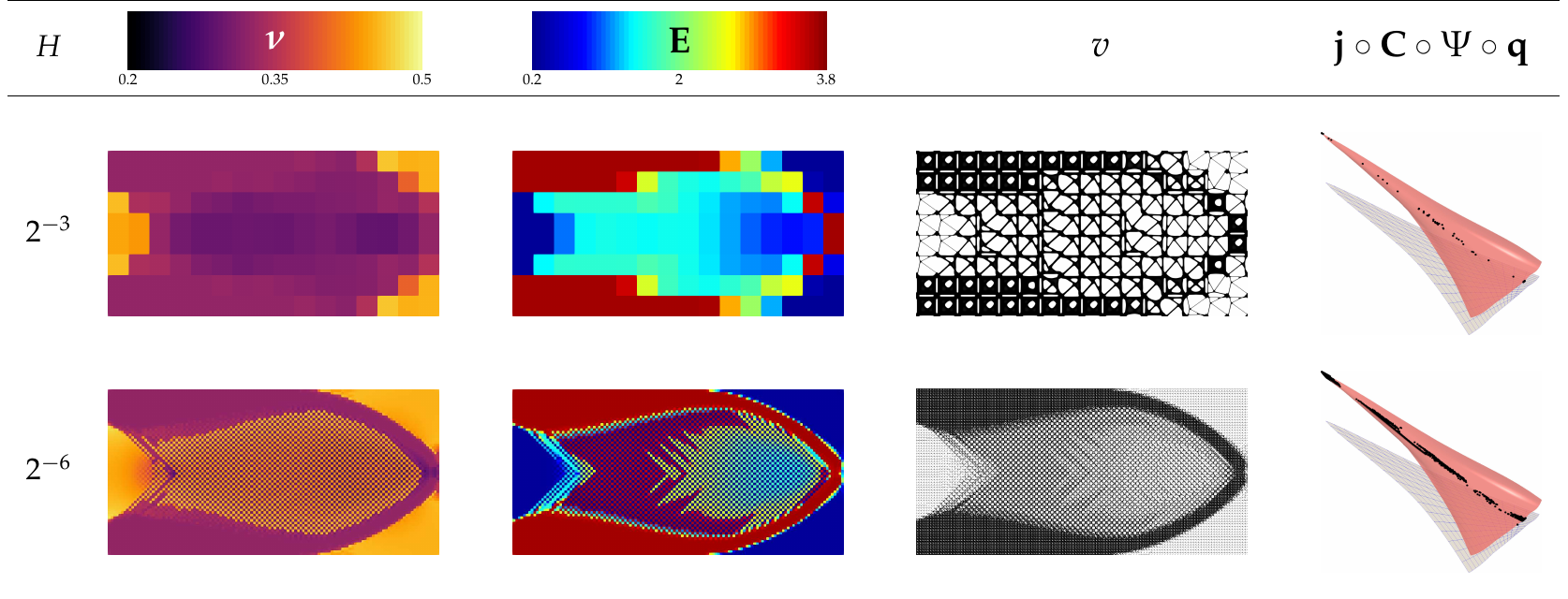}} 
\caption[2D Cantilever, bridges corner, vol 1]{
Macroscopic and microscopic solutions are visualized for the cantilever problem and a volume constraint $\Volume_H = 1$ and the admissible fine-scale structures with material bridges at the corners.}
\label{fig:TwoScaleCantilverBridgesCornerVol1}
\end{figure}

\begin{figure}[htbp]
 {\includegraphics[width=\textwidth]{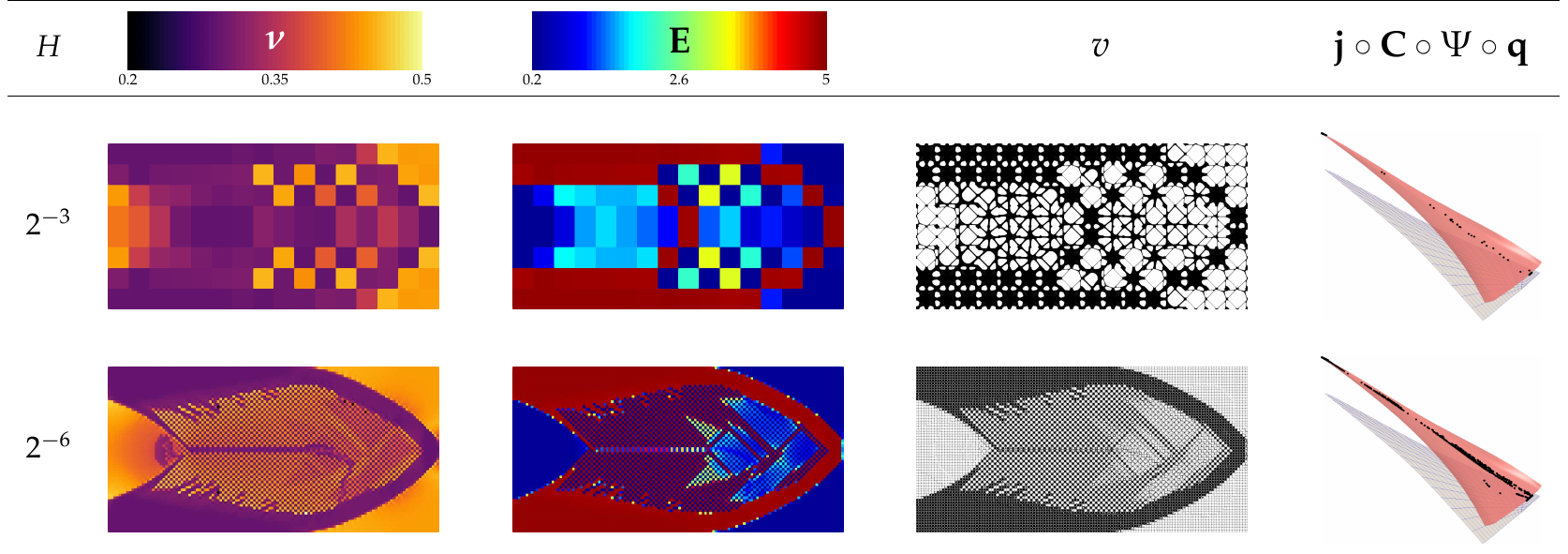}} 
\caption[2D Cantilever, bridges mid and corner, vol 1]{
Macroscopic and microscopic solutions are rendered for the cantilever problem for a volume constraint $\Volume_H = 1$
and the admissible fine-scale structures with material bridges both at the corners and midfaces.
}
\label{fig:TwoScaleCantilverBridgesMiddleAndCornersVol1}
\end{figure}

\begin{figure}[htbp]
\centering{
 {\includegraphics[width=\textwidth]{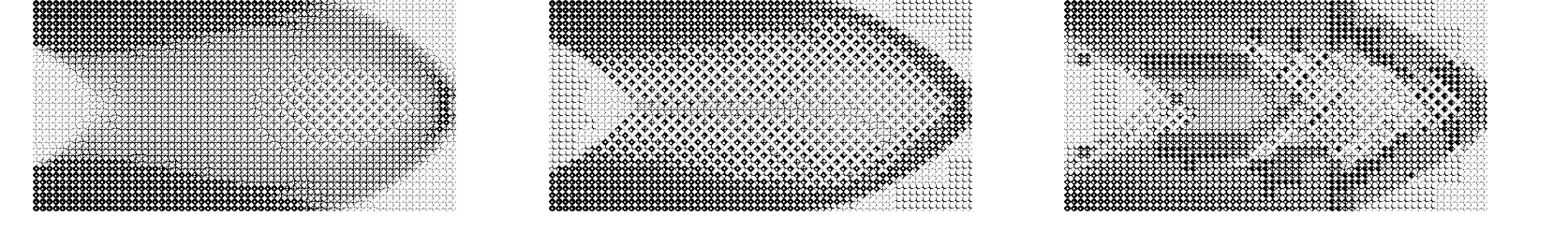}} 
 \begin{tabular}{|c | c |}
 \hline
  $\widehat{c_P}$ & $0$ \\
  \hline 
  $\mathcal{J}_{\text{compl}}$ &  0.098566 \\
  $L^{\sigma}$                &  759.374 \\  
  \hline
 \end{tabular}
 \hspace*{0.125\textwidth}
  \begin{tabular}{| c | c |}
 \hline
  $\widehat{c_P}$ & $10^{-4}$ \\
  \hline 
  $\mathcal{J}_{\text{compl}}$ &  0.103957896 \\
  $L^{\sigma}$                 &  669.9991    \\  
  \hline
 \end{tabular}
  \hspace*{0.125\textwidth}
  \begin{tabular}{|c | c |}
 \hline
  $\widehat{c_P}$ & $10^{-3}$ \\
  \hline 
  $\mathcal{J}_{\text{compl}}$  & 0.12192695 \\
  $L^{\sigma}$                  & 638.023886 \\  
  \hline
 \end{tabular}
 }
\caption[2D Cantilever, bridges mid, vol 0.75, compare cp]{
Optimal two-scale structures are depicted for the cantilever problem with admissible fine-scale structures based on material 
bridges at midfaces, a volume constraint $\Volume_H = 0.75$, and a macroscopic grid size $\macrogridsize = 2^{-5}$. 
From left to right we have chosen $\widehat{c_P} = 0, \, 10^{-4}, \, 10^{-3}$ to penalize the interface for the total cost functional.
Whereas, for increasing values of $\widehat{c_P}$, the interfacial energy of the corresponding optimal structure is decreasing, the compliance and the total macroscopic energy is increasing.
}
\label{fig:TwoScaleCantilverBridgesMiddleVol075Compare}
\end{figure}

\paragraph*{Minimization of the potential energy in 2d.}
First, we consider the potential energy as objective functional.
In the following, let $\domain = [0,2] \times [0,1]$ be the computational domain for the macroscopic problem.
The Dirichlet boundary is given by $\Gamma_D = {0} \times [0,1]$ and we apply a force $f = (0, -10 \, \chi_{[1.95,2] \times [0.45,0.55]} )^T$.
Here, we always fix a certain amount of hard material, which we denote by $\Volume_H$.
More precisely, by choosing $\widehat{c_P} = 0$, we consider the optimization problem
\begin{align}
 \argmin_{\mathbf{q} \ : \ \int_\domain \mathbf{j}_{\text{ref}}(\mathbf{q}(x)) \d x = \Volume_H}
 \mathcal{J}_{\text{compl}}[\displaceMacro[\Etensor(\Psi \circ \mathbf{q})]]
\end{align}
In Figure~\ref{fig:TwoScaleCantilverBridgesMiddleVol1}, we depict the results for a volume constraint $\Volume_H = 1$
and the bridge set $B^1$ given in the middle of the edges.
In particular, for macroscopic grid sizes $\macrogridsize=\epsilon=2^{-3},2^{-4},2^{-5},2^{-6}$, we show full two-scale structures.
Moreover, we explore these results for finer grid sizes $\macrogridsize=\epsilon=2^{-7},2^{-8}$.
In Figure~\ref{fig:TwoScaleCantilverBridgesCornerVol1} and Figure~\ref{fig:TwoScaleCantilverBridgesMiddleAndCornersVol1}, we compare the results for the different bridge sets.
Remember from the numerical results on the micro-scale (\cf Figure~\ref{fig:CellOptVol}) that for the other bridge sets the ranges of the effective Young's modulus are essentially higher.
These ranges are of course exhausted in the macroscopic optimization, in particular note these different ranges for $E$ in the described figures.
Furthermore, in Figure~\ref{fig:TwoScaleCantilverBridgesMiddleVol075Compare}, we study the effect of the interface penalization parameter $\widehat{c_P}$ for a fixed volume constraint $\Volume_H=0.75$. 

\begin{figure}[htbp]
 {\includegraphics[width=\textwidth]{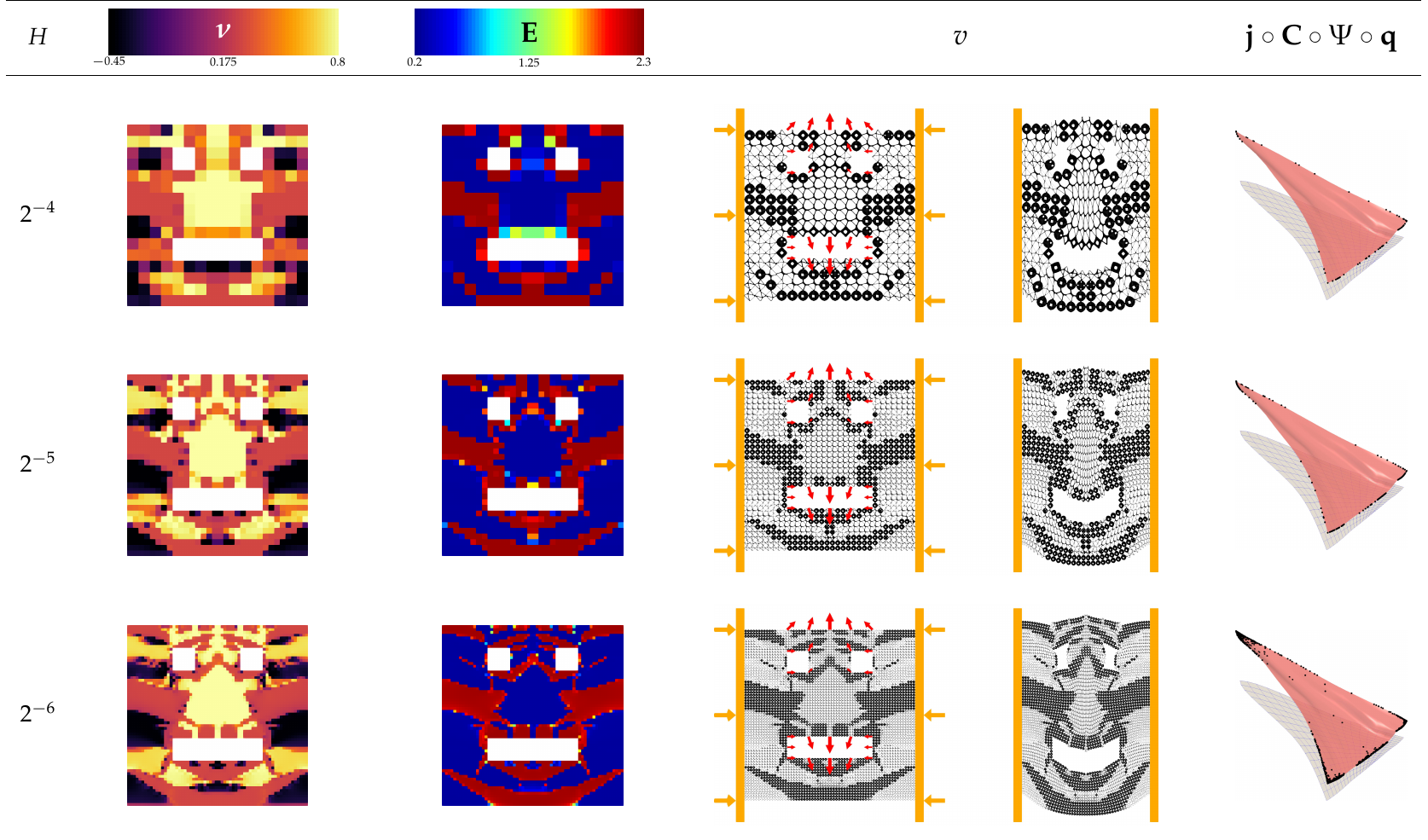}} 
\caption[2D Smiley]{Shape optimization solution for a tracking type functional enforcing a smiley in 2D under uniform horizontal displacements on the left and right boundary. The tracking domain is the union of all cells with a red arrow indicating the displacement $U^0$.}
\label{fig:Smiley2D}
\end{figure}

\begin{figure}[htbp]
 {\includegraphics[width=\textwidth]{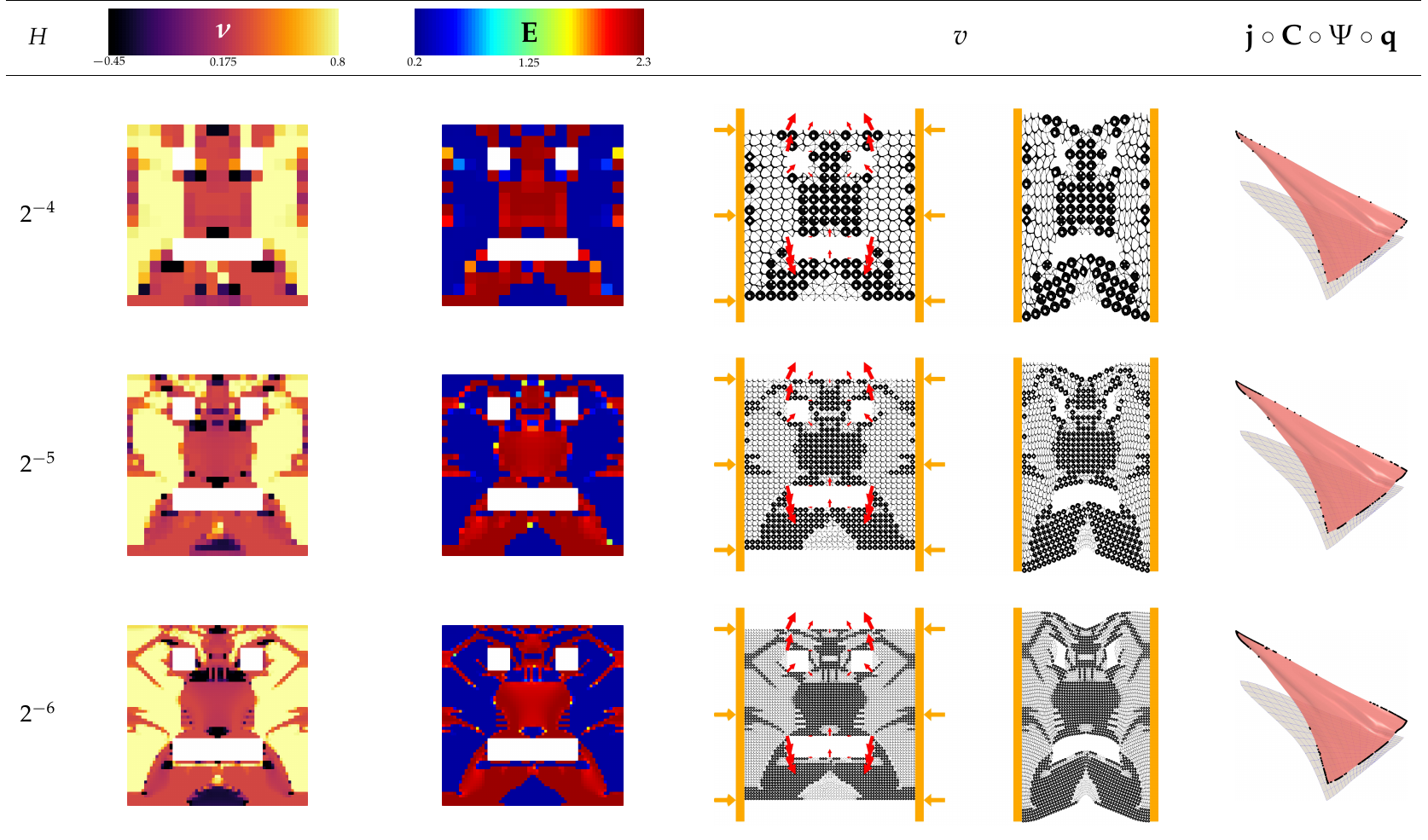}}
\caption[2D Frowney]{Shape optimization solution for a tracking type functional (\cf~Figure~\ref{fig:Smiley2D}) enforcing a frowney in 2D under uniform horizontal displacements on the left and right boundary.}
\label{fig:Frowney2D}
\end{figure}

\paragraph*{Tracking type functional in 2d.}
Now, we choose the domain  
\begin{align} \label{eq:domainSmiley2D}
 \domain 
 = \domain_{2D} 
 = [0,1]^2 \setminus \left(
   (\tfrac{1}{4},\tfrac{3}{4}) \times (\tfrac{1}{4},\tfrac{3}{8}) \cup  (\tfrac{1}{4},\tfrac{3}{8}) \times (\tfrac{3}{4},\tfrac{7}{8}) \cup (\tfrac{5}{8},\tfrac{3}{4}) \times (\tfrac{3}{4},\tfrac{7}{8})
   \right) \, ,
\end{align}
where we apply a prescribed compression displacement
\begin{align}
 \displaceMacro^\partial_1(0,x_2) = \tfrac{1}{8}, \quad \displaceMacro^\partial_1(1,x_2) = -\tfrac{1}{8} \, .
\end{align}
To obtain a unique solution of the elastic problem, we fix the integral $\int_\domain \displaceMacro_2 \d x = 0$.
In Figure~\ref{fig:Smiley2D}, on a localized tracking domain, we show the result for tracking type displacement to generate a smiley in the deformed configuration.
Similarly, in Figure~\ref{fig:Frowney2D}, we show the result to generate a frowney.
Note that, compared to the results for the Cantilever beam, where the point cloud on the admissible surface is essentially supported on a line through the middle,
here, the range of Poisson's ratio is fully exhausted in both directions. 

\begin{figure}[htbp]
 {\includegraphics[width=\textwidth]{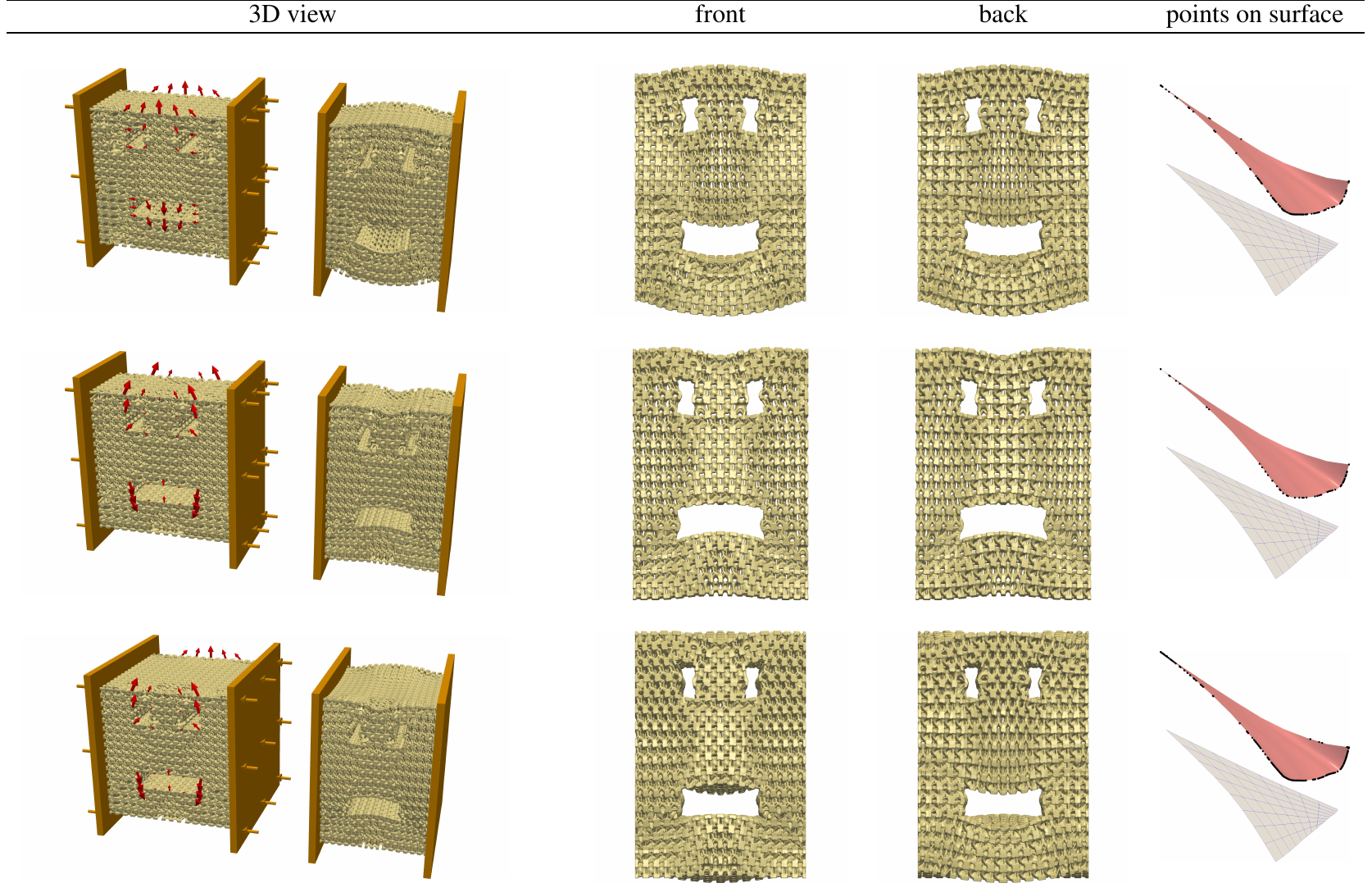}}
\caption[3D Smiley]{We show the shape optimization solution for several tracking type functionals in $3D$. 
From top to bottom, we depict a smiley, a frowney, and an object which is angry at the front face and laughing at the back face.
Here, the macroscopic grid size is always given by $\macrogridsize = 2^{-4}$.}
\label{fig:Smiley3DAllRef1}
\end{figure}

\paragraph*{Tracking type functional in 3d.}
Finally, we consider similar tracking type problems in 3D.
Given the domain $\domain_{2D}$ as defined in \eqref{eq:domainSmiley2D} and $Z>0$, we set $\domain = \domain_{2D} \times [0,Z]$.
Then we apply a compression displacement
\begin{align}
 \displaceMacro^\partial_1(0,x_2,x_3) = \tfrac{1}{8}, \quad \displaceMacro^\partial_1(1,x_2,x_3) = -\tfrac{1}{8} 
\end{align}
which does not constrain in-plane deformation on the $(x_2,x_3)$--boundary planes.
Thus, in addition to the integral constraints $\int_\domain \displaceMacro_2 \d x = \int_\domain \displaceMacro_3 \d x = 0$, we fix the rotation in these boundary planes by enforcing $\int_\domain (\curl \displaceMacro)_1 \d x = 0$ to obtain a unique solution of the elastic problem.
In analogy to the two-dimensional case, we show in Figure~\ref{fig:Smiley3DAllRef1} the resulting optimized shape for two different tracking type functionals, which enforce either a smiley or a frowney, 
in both examples we use $Z=\tfrac{1}{2}$.
Furthermore, we consider for $Z=1$ a tracking type problem generating a frowney on the front face and a smiley on the back face.

\begin{figure}[htbp]
\centering{
 \begin{minipage}{0.45\textwidth}
    \includegraphics[height=4.5cm]{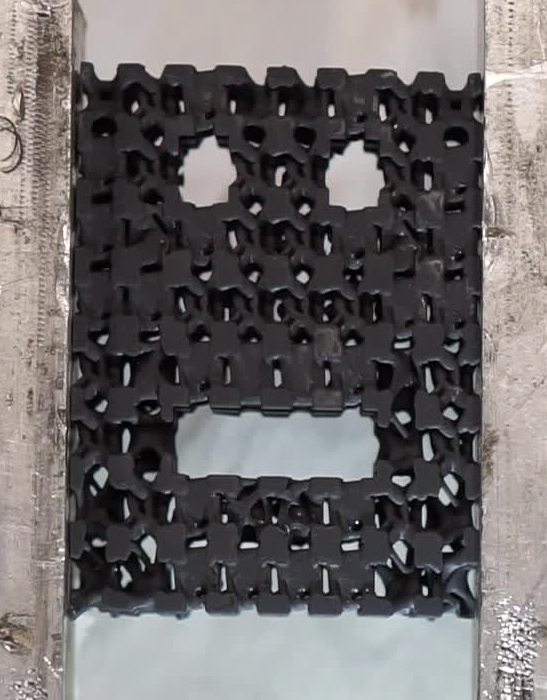}
    \hspace*{0.2cm}
    \includegraphics[height=4.5cm]{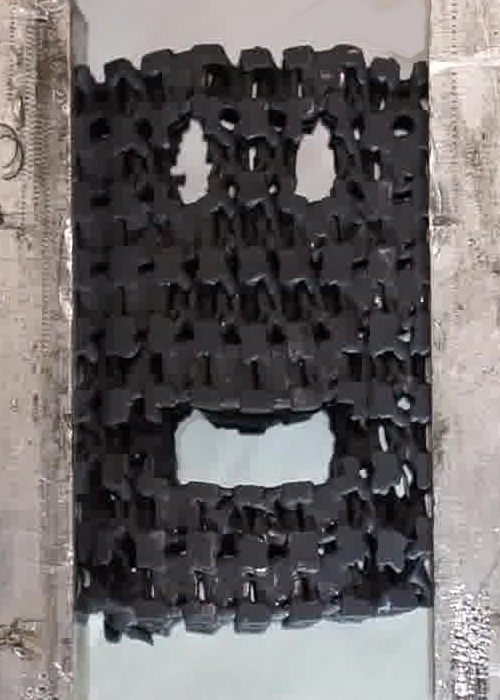}
  \end{minipage}
  \hspace*{0.02\textwidth}
 \begin{minipage}{0.45\textwidth}
    \includegraphics[height=4.5cm]{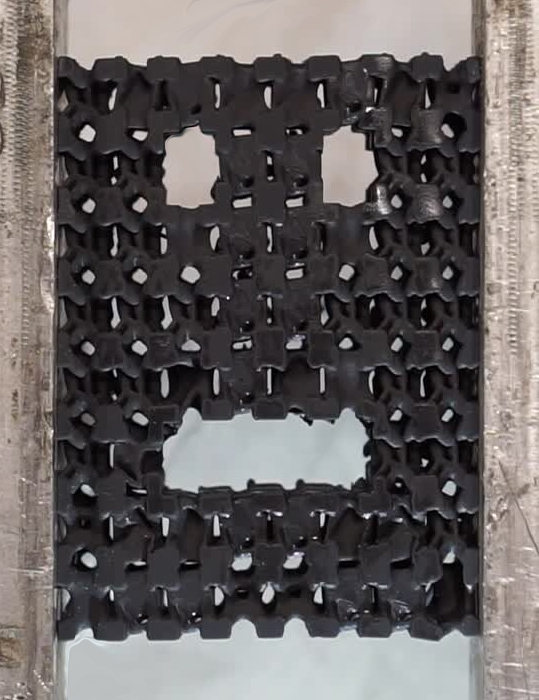}
    \hspace*{0.2cm}
    \includegraphics[height=4.5cm]{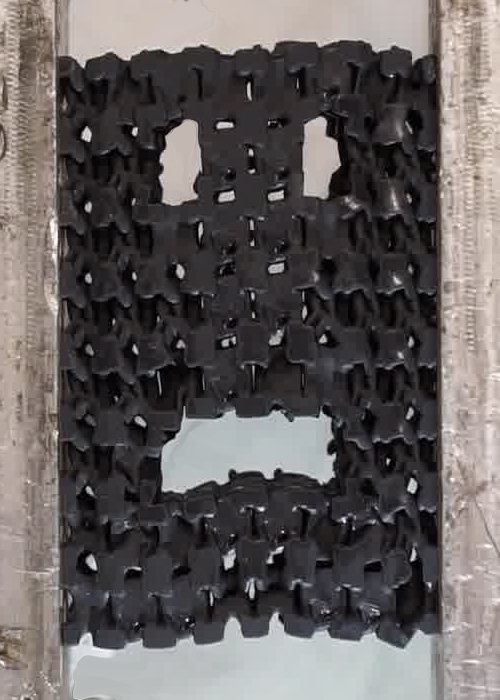}
  \end{minipage}
  }
\caption[3D Print]{
Photos of two $3D$ printed objects resulting from the free material optimization for the two tracking type functionals from 
Fig.~\ref{fig:Smiley3DAllRef1}:
On the left, the undeformed and deformed configuration leading to a smiley is depicted.
On the right, the same is shown for the tracking type functional tailored to enforce a frowney.
}
\label{fig:Smiley3DPrint}
\end{figure}

As a proof of concept for the actual manufacturability, we printed the computed workpieces 
with an SLA printer with relatively soft, black resin. 
To this end, we considered the shape optimization for the tracking type functionals above generating a smiley and a frowney, respectively.
Due to limitations of the 3D printer we considered 
$Z=\tfrac{1}{4}$ and a larger macroscopic grid size. 
The one-dimensional compression in horizontal direction is obtained in a vise with sufficiently large 
vise jaws covering the left and the right side of the object.
The total size of both objects is  $8 \times 8 \times 2 \ cm^3$,
the volume of a single cell is  $1 \ cm^3$.
Figure~\ref{fig:Smiley3DPrint} depicts photos of the undeformed and the deformed configurations of the objects. 

\section{Conclusion} \label{sec:conclusion}
We have discussed a new two-scale approach for elastic shape optimization of fine-scale structures in additive manufacturing.
In an offline phase, a database of admissible cell microstructures is computed, 
where the manufacturability constraint is implemented by prescribing common material bridge sets between neighbouring cells.
For a given effective (isotropic) elasticity tensor, these structures are optimized  
\wrt a cost functional combining the volume of the printed material and the length (or the area) of the interfaces
within a single cell.
The online phase consists in an efficient
restricted
free material optimization on the macroscopic scale by selecting cell-wise 
the optimal material properties within the realizable set that was identified in the offline phase. 
This set is parametrized in the space of elasticity constants using a spline-based sampling strategy.
As possible macroscopic cost functionals a compliance energy and a tracking-type functional are discussed.
Isotropic effective elasticity tensors are chosen to reduce the computations to a two parameter sampling in the 
offline database. 
However, the proposed method to recover a specific effective elasticity tensor with optimal cost could be used for more general elasticity tensors as well.
Here, a hierachical B-spline interpolation on sparse grids as proposed in \cite{VaHuSt20} might help to cope with higher dimensional admissible sets of elasticity tensors.

As we have mentioned in Section~\ref{sec:intro}, many two-scale approaches have been applied recently to generate manufacturable workpieces.
Note that our approach takes into account already optimized microstructures for fixed isotropic effective elasticity tensors.
In contrast, Panetta~\etal~\cite{PaZhMa15} used an isotropic filtering on explicitly constructed truss geometries.
With the intention to approximate nested laminates, Allaire~\etal~\cite{AlGePa19,GeAlPa20} considered microstructures with rectangular holes.
Full shape optimization of the microstructures has been performed by Ferrer~\etal~\cite{FeCaHe18} \wrt stresses.
Similar to our approach, Schury~\etal~\cite{ScStWe12} and Zhu~\etal~\cite{ZhSkCh17} tried to recover given effective elasticity tensors,
where both used tracking type functionals.
Instead, we have demonstrated that our optimization scheme allows a precise recovery of a large point cloud of isotropic effective elasticity tensors.
Furthermore, manufacturability is directly guaranteed due to the predefined material bridges, which allow a proper transfer of stresses.

\section*{Acknowledgments}
We thank Patrick Dondl for useful discussions, in particular in connection with 3D printing.
We acknowledge support by the Computing Center at the University of Bonn with respect to the 3D printing.
This work was partially supported
by the Deutsche Forschungsgemeinschaft through project 211504053/SFB1060.

\bibliographystyle{siamplain}

\def\polhk#1{\setbox0=\hbox{#1}{\ooalign{\hidewidth
  \lower1.5ex\hbox{`}\hidewidth\crcr\unhbox0}}}

\end{document}